\documentclass[12pt]{amsart}
\usepackage{amscd}
\usepackage{graphics}
\usepackage{afterpage}
\input prepictex   \input pictex    \input postpictex

\def\bulletsymb{\large\bf.}

\def \Ext{\operatorname{Ext}}
\def \End{\operatorname{End}}
\def \Im{\operatorname{Im}}
\def \op{{\rm op}}
\def \Rad{\operatorname{Rad}}
\def \Hom{\operatorname{Hom}}
\def \mod{\operatorname{mod}}
\DeclareMathOperator\sub{sub}
\def \fac{\operatorname{fac}}

\def \rad{\operatorname{rad}}
\def \soc{\operatorname{soc}}
\def \lto#1{\;\mathop{\longrightarrow}\limits^{#1}\;} 
 
\def \sto#1{\;\mathop{\to}\limits^{#1}\;}
\def \arr#1#2{\arrow <2mm> [0.25,0.75] from #1 to #2}
\def \prep#1#2#3#4#5#6#7#8{\vcenter{\beginpicture
    \setcoordinatesystem units <5mm,5mm>
    \put{$\scriptstyle #1$} at 1 1
    \put{$\scriptstyle #2$} at 0 0
    \put{$\scriptstyle #3$} at 2 0
    \put{$\scriptstyle #4$} at 1 -1
    \arr{.3 .3}{.7 .7}
    \arr{1.7 .3}{1.3 .7}
    \arr{.7 -.7}{.3 -.3}
    \arr{1.3 -.7}{1.7 -.3}
    \put{${\scriptstyle #5}\;$} [br] at .4 .6
    \put{$\;{\scriptstyle #6}$} [bl] at 1.6 .6
    \put{${\scriptstyle #7}\;$} [tr] at .4 -.6
    \put{$\;{\scriptstyle #8}$} [tl] at 1.6 -.6
    \endpicture}}
\def \pmap#1#2#3#4{\scriptstyle #2{{\scriptstyle #1}\atop{\scriptstyle#4}}#3}

  \def\eyesm{\put{} at 0 -.1 \multiput{{\Large\rm .}} at -.1 .1  .1 .1 / }
  \def\eyesh{\put{} at 0 -.35 \multiput{{\Large\rm .}} at -.1 .35  .1 .35 / }
  \def\eyesl{\put{} at 0 .25 \multiput{{\Large\rm .}} at -.1 -.15  .1 -.15 / }
  \def\nosem{\plot 0 .1  0 -.1 / }
  \def\noseh{\put{} at 0 -.35 \plot 0 .35  0 .15 / }
  \def\nosel{\put{} at 0 .35 \plot 0 -.15  0 -.35 / }
  \def\mouthm{\put{} at 0 .2 \circulararc 60 degrees from 0 -.2 center at 0 0
    \circulararc -60 degrees from 0 -.2 center at 0 0 }
  \def\mouthh{\put{} at 0 -.2 \circulararc 60 degrees from 0 .05 center at 0 .25
    \circulararc -60 degrees from 0 .05 center at 0 .25 }
  \def\mouthl{\put{} at 0 .4 
    \circulararc 60 degrees from 0 -.45 center at 0 -.25
    \circulararc -60 degrees from 0 -.45 center at 0 -.25 }
  
\def\sq{\plot 0 0  1 0  1 1  0 1  0 0 /}
\def\oneo{\put{$\smallsq5$} at -.25 0 }
\def\oo{\multiput{$\smallsq5$} at -.25 -.25  -.25 .25 / }
\def\smallsq#1{\plot 0 0  0.#1 0  0.#1 0.#1  0 0.#1  0 0 /}
  
\def\qed{\phantom{m.} $\!\!\!\!\!\!\!\!$\nolinebreak\hfill\checkmark}
\newcommand\Zmod[1]{\mathbb Z/p^#1}

\newtheoremstyle{mytheorems}{9pt}{6pt}{\itshape}{0pt}{\sc}{.}{ }{}
\newtheoremstyle{myremarks}{6pt}{3pt}{\normalfont}{0pt}{\it}{.}{ }{}

\theoremstyle{mytheorems}
\newtheorem{theorem}{Theorem}[section]
\newtheorem{lemma}[theorem]{Lemma}
\newtheorem{corollary}[theorem]{Corollary}
\newtheorem{proposition}[theorem]{Proposition}

\theoremstyle{myremarks}
\newtheorem*{notation}{Notation}
\newtheorem{example}[theorem]{Example}
\newtheorem{observation}[theorem]{Observation}
\newtheorem{remark}[theorem]{Remark}
\newtheorem*{definition}{Definition}

\begin{document}
                                %
\noindent{\phantom{\footnotesize [s-finite1k, September 12, 2006]}}
                                %
                                %
\vglue1truecm
\centerline{\Large Systems of Submodules}
\medskip\centerline{\Large and a Remark by M.~C.~R.~Butler}
\bigskip
\centerline{By}
\medskip
\centerline{\sc Markus Schmidmeier}

\bigskip\medskip

\centerline{\parbox{10cm}{\footnotesize 
    {\it Abstract.} 
    Fix a poset $\mathcal P$ and a natural number $n$.
    For various commutative local rings $\Lambda$, each of 
    Loewy length $n$, consider the category 
    $\textrm{sub}_\Lambda\mathcal P$ of $\Lambda$-linear submodule
    representations of $\mathcal P$.
    We give a criterion for when the underlying translation quiver
    of a connected component of the Auslander-Reiten
    quiver of $\sub_\Lambda\mathcal P$ is independent of the choice
    of the base ring $\Lambda$.
    If $\mathcal P$ is the one-point poset and 
    $\Lambda=\mathbb Z/p^n$, then $\textrm{sub}_\Lambda\mathcal P$
    consists of all pairs $(B;A)$ where $B$ is a finite abelian $p^n$-bounded group 
    and $A\subset B$ a subgroup.
    We can respond to a remark by M.~C.~R. Butler
    concerning the first occurence of parametrized 
    families of such subgroup embeddings.}}

\renewcommand{\thefootnote}{}

\section{Introduction.}
Let $\Lambda$ be a commutative (local artinian) uniserial ring, for example
$\Lambda = k[x]/x^n$ where $k$ is a field or $\Lambda=\mathbb Z/p^n$ where
$p$ is a prime number, and let $\mathcal P$ be a finite partially ordered set
(poset).
We consider systems $(M_*,(M_i)_{i\in\mathcal P})$ where $M_*$ is a finitely
generated $\Lambda$-module and each $M_i\subset M_*$ is a submodule such that
$M_i\subset M_j$ holds for $i<j$ in $\mathcal P$.
The systems $(M_*,(M_i)_{i\in\mathcal P})$ form the objects of a category
$\sub_\Lambda\mathcal P$; here homomorphisms from $(M_*,(M_i))$
to $(N_*,(N_i))$ are given by $\Lambda$-linear maps $f:M_*\to N_*$ which
satisfy $f(M_i)\subset N_i$ for each $i\in\mathcal P$. 
This defines an additive category which has the Krull-Remak-Schmidt property,
so any object has a unique decomposition into indecomposable ones.
Thus, the focus is on the classification of the indecomposable objects and 
a description of the homomorphisms between them.

\smallskip
Given two commutative uniserial rings $\Lambda$ and $\Delta$ of the same
Loewy length $n$, then one may ask how the categories $\sub_\Lambda\mathcal P$
and $\sub_\Delta\mathcal P$ are related: Do they have the same number of 
isomorphism classes of indecomposable objects? If so, is there a bijection
between the isomorphism classes which preserves 
(A)~combinatorial data that describe the objects and
(B)~combinatorial data which describe the homomorphisms between them?
We will describe objects in terms of their {\it type:} 
For a $\Lambda$-module $M_*=\Lambda/p^{m_1}\oplus\cdots\Lambda/p^{m_t}$, 
$p$ a radical generator for $\Lambda$, the type is the partition 
$t(M_*)=(m_1,\ldots,m_t)$; then for a system $M=(M_*,(M_i))$, the type
is given by the tuple $t(M)=(t(M_*),t(M_i))$. 
For the description of homomorphisms we consider the combinatorial data given by the
Auslander-Reiten quiver of the category. 

\smallskip
The categories of type $\sub_\Lambda\mathcal P$ studied for example in \cite{richman}
and \cite{rump} have finite representation type; it is shown there that 
the classification of the 
indecomposable objects does not depend on the choice of the base ring $\Lambda$.
On the other hand covering theory, which is used for example in
\cite{rs-inv} and \cite{ms-bounded} for certain categories of finite and tame 
infinite representation type
delivers both the full list of the indecomposable modules and the structure of
the Auslander-Reiten quiver.
However, this method requires that the uniserial base ring $\Lambda$ is an algebra over
a field, hence of type $\Lambda=k[x]/x^n$.

\smallskip
If $\Delta$ is another commutative uniserial ring of length $n$, then we can use
the main result in this manuscript as a criterion for the existence of a bijection
between the indecomposable objects in $\sub_\Lambda\mathcal P$ and $\sub_\Delta\mathcal P$
which preserves (A) and (B).

\begin{theorem} \label{introthm}
  Suppose $\Lambda,\Delta$ are commutative
  uniserial rings of the same length $n$, $\mathcal P$ is a finite poset, 
  $\Gamma_\Lambda$ and $\Gamma_\Delta$ are
  connected components 
  of the Auslander-Reiten quivers of $\sub_\Lambda\mathcal P$ and
  $\sub_\Delta\mathcal P$, and $\mathcal L_\Lambda$ and $\mathcal L_\Delta$
  are slices in $\Gamma_\Lambda$ and $\Gamma_\Delta$, respectively, such that the following 
  conditions are satisfied:
  \begin{enumerate}
  \item The slices $\mathcal L_\Lambda$ and $\mathcal L_\Delta$ are isomorphic as 
    valued graphs.
  \item Points in $\mathcal L_\Lambda$ and $\mathcal L_\Delta$ which correspond to each other
    under this isomorphism represent indecomposable objects of the same type.
  \item Each indecomposable object represented by a point in $\mathcal L_\Lambda$ or in
    $\mathcal L_\Delta$ is determined uniquely, up to isomorphism, by its type.
  \end{enumerate}
  Then the components $\Gamma_\Lambda$ and $\Gamma_\Delta$ are isomorphic
  as valued translation quivers, 
  and points which correspond to each other under this isomorphism
  represent indecomposable objects of the same type.
\end{theorem}

\smallskip
{\it Subgroups of abelian groups.}
Our research was motivated by a remark by M.C.R.\ Butler after the authors talk
on categories of type $\mathcal S_m(k[x]/x^n)$ at the ICRA XI in 2004 in Patzcuaro,
Mexico. 
For $\Lambda$ a commutative uniserial ring of length $n$ and $m\leq n$, the 
category $\mathcal S_m(\Lambda)$ consists of all pairs $(A\subset B)$ where
$B$ is a (finitely generated) $\Lambda$-module and $A$ a submodule of $B$
which is annihilated by $\rad^m\Lambda$.
Thus $\mathcal S_n(\Lambda)=\sub_\Lambda\mathcal P$ for $\mathcal P$ the one point
poset, and for each $m<n$, $\mathcal S_m(\Lambda)$ is a full subcategory.
Using coverings, it has been shown in \cite{ms-bounded} for which pairs $(m,n)$ 
the category $\mathcal S_m(\Lambda)$
where $\Lambda=k[x]/x^n$  has finite, tame or wild representation type.  
It turns out that the two pairs where $(m,n)=(4,6)$ or $(3,7)$ 
mark the first occurances of parametrized families of indecomposable objects. 

\smallskip
Since Birkhoff \cite{birkhoff} the sub{\it group} categories $\mathcal S_m(\Lambda)$ 
where $\Lambda=\mathbb Z/p^n$ for $p$ a prime number have attracted a lot of interest.
Here we are dealing with the possible embeddings of $p^m$-bounded subgroups in
$p^n$-bounded finite abelian groups.
It is of particular interest, as Butler points out, to detect the first occurances of 
parametrized families of indecomposable objects in the case of subgroup embeddings. 
It is shown in \cite{birkhoff} and in \cite{ms-bounded} that 
parametrized families of indecomposable objects do occur 
in $\mathcal S_4(\mathbb Z/p^6)$ and in $\mathcal S_3(\mathbb Z/p^7)$, respectively.
But are these the {\it first} occurances 
(as in the case of modules over the polynomial ring)?  
In the group case we know from \cite{richman} 
that the category $\mathcal S_5(\mathbb Z/p^5)$ has finite
type, so there are no parametrized families of indecomposables in any of the 
categories $\mathcal S_m(\mathbb Z/p^n)$ for $m\leq n\leq 5$.  
Also if $m\leq2$ then it is known in group theory \cite{bhw} that for any $n$,
each indecomposable object $(A\subset B)$ in $\mathcal S_m(\mathbb Z/p^n)$ 
is determined uniquely by its type
--- again there are no families.  
So there remains a single case to decide, namely where $m=3$ and $n=6$,
and this will be our main example.

\bigskip
{\it Related results.}
For previous results in the representation theory of a $\Lambda$-category
where the base ring $\Lambda$ is not necessarily a field or a finite dimensional algebra,
see in particular \cite{simson} for the case where $\mathcal P$ is a
chain, \cite{avino-bautista} for modules over the group ring 
$\mathbb Z/p^n\; C_p$, 
and \cite{rump} for lattices over tiled orders.
In \cite{ps} it is shown that for suitable categorical ideals $\mathcal I$ and $\mathcal J$,
the categories $\mathcal S_2(k[x]/x^n)/\mathcal I$ and $\mathcal S_2(\mathbb Z/p^n)/\mathcal J$
are equivalent categories.  Moreover, according to \cite{rs-wild},
the category $\mathcal S_m(\Lambda)$ 
is controlled wild whenever the base ring $\Lambda$ has
length at least 7 and $m\geq 4$.

\bigskip
{\it The contents of each section.}
In order to reconstruct the underlying (valued) translation quiver of an
Auslander-Reiten quiver from a given slice, it is essential that 
projective objects, injective objects, the summands of the first term in a sink map 
for a projective object, and the summands of the 
last term of a source map for an injective object can be detected by their types.
This is the case, as we see in the next section if $\Lambda$ is a commutative
uniserial ring and we are dealing with a category of type $\sub_\Lambda\mathcal P$
or $\mathcal S_m(\Lambda)$. 

\smallskip
Given a short split exact sequence, then clearly, the type of the middle term is the 
union of the types of the end terms.  This is also true for almost split sequences,
with finitely many exceptions.  In Section~\ref{section-nonsplit} 
we describe these exceptional sequences. It turns out that the starting and end terms 
of these sequences are indecomposable objects determined uniquely by their types. 

\smallskip
As a consequence if $\Gamma$ is a connected component of the Auslander-Reiten quiver
and $\mathcal L$ a slice in $\Gamma$, then the type of each module in $\Gamma$ is 
determined uniquely by the types of the modules in $\mathcal L$. 
This implies, as we show in Section~\ref{section-isomorphy} that whenever
the conditions in Theorem~\ref{introthm} apply, the combinatorial structure of the
component of the 
Auslander-Reiten quiver does not depend on the base ring chosen. 
For illustration we consider the category of all chains of length 3 of $\Lambda$-modules
where $\Lambda=k[x]/x^2$ or $\Lambda=\mathbb Z/p^2$. This case 
is particularly straightforward
as the slice is obtained from the sequences studied in Section~\ref{section-nonsplit}.

\smallskip
In general it may require a lot of legwork to verify that a particular short exact sequence
is almost split.  A test for Auslander-Reiten sequences 
from \cite{ars} applies also to our situation, as we show in
Section~\ref{section-test}.  With the help of this test we obtain the slice needed
in Theorem~\ref{introthm} to 
demonstrate that the structure of the Auslander-Reiten quiver for $\mathcal S_3(\Lambda)$ 
is independent of the choice of the commutative uniserial ring $\Lambda$ of length~6.

\smallskip 
{\it Notation and Remarks.} 
For terminology related to Auslander-Reiten sequences, relative homological algebra and
coverings we refer the reader to \cite{ars}, \cite{as}, and \cite{bongartz-gabriel}.
Some of the results in this manuscript have been presented at the 
2005 Oberwolfach meeting on representation theory of finite dimensional
algebras~\cite{oberwolfach},
and at the 2005 regional meeting of the AMS in Santa Barbara.

\section{Projectives and injectives.}\label{projinj}

We determine the indecomposable projective and injective objects in the submodule category
$\sub_\Lambda\mathcal P$ and their respective sink and source maps.
It turns out that whenever $\Lambda$ is a commutative uniserial ring, 
then the first term and the
last term of each such map is an indecomposable module, which is determined uniquely
by its type, or zero.

  \smallskip
  For $\mathcal P$ a finite poset, 
  denote by ${\mathcal P}^*$ and by ${\mathcal P}^0$ the 
  poset obtained from $\mathcal P$ 
  by adding a largest and a smallest element, respectively.
  The incidence algebra $\Lambda\mathcal P^*$ of the poset $\mathcal P^*$ has free $\Lambda$-basis
  $\{q_{ij}|i,j\in\mathcal P^*, i\leq j\}$ and multiplication is given by the formula
  $q_{ij}q_{k\ell}=\delta_{jk}q_{i\ell}$.
  We identify $\sub_\Lambda\mathcal P$ with the full subcategory of $\mod\Lambda\mathcal P^*$
  of all (right) modules $M$ such that for each pair $i,j\in\mathcal P^*,i\leq j,$ the 
  multiplication map $Mq_{ii}\to Mq_{jj}, m\mapsto mq_{ij},$ is monic.
  Each indecomposable projective $\Lambda\mathcal P^*$-module has the form $q_{ii}\Lambda\mathcal P^*$
  for some $i\in\mathcal P^*$ and hence the module itself and its radical both lie in the
  subcategory $\sub_\Lambda\mathcal P$.

  \smallskip
  Similarly we identify $\fac_\Lambda\mathcal P$, the category of all finitely generated
  $\Lambda$-linear factor space representations of $\mathcal P$, with the full subcategory
  of $\mod\Lambda\mathcal P^0$ of all modules for which each of the maps
  $Mq_{ii}\to Mq_{jj}, m\mapsto mq_{ij},$ where $i,j\in\mathcal P^0$, $i\leq j$, is onto.
  The self duality $\Hom_\Lambda(-,I):\mod\Lambda\to \mod\Lambda$ given by the injective 
  envelope $I=E(\Lambda/\rad\Lambda)$ \cite[30.6]{af} gives rise to a duality 
  $\mod\Lambda(\mathcal P^0)^\op\to\mod\Lambda\mathcal P^0$. 
  Hence the injective $\Lambda\mathcal P^0$-modules are under control; they and their socle
  factors all lie in the subcategory $\fac_\Lambda\mathcal P$.

  \medskip
  Neither $\sub_\Lambda\mathcal P$ nor $\fac_\Lambda\mathcal P$ is an abelian category,
  unless $\mathcal P$ is the empty poset.  However, both $\sub_\Lambda\mathcal P$ and
  $\fac_\Lambda\mathcal P$ are full exact subcategories of the module categories
  $\mod \Lambda\mathcal P^*$ and $\mod\Lambda\mathcal P^0$, respectively.

\begin{lemma}
  \begin{enumerate}
    \smallskip \item[1a)] Every projective 
    $\Lambda{\mathcal P}^*$-module
    is a projective subspace representation 
    of $\mathcal P$.
    
  \item[1b)] Every injective $\Lambda{\mathcal P}^0$-module 
    is an injective
    factorspace representation of $\mathcal P$.
    
    \smallskip
  \item[2a)] For each indecomposable projective object 
    $P$ in $\sub_\Lambda{\mathcal P}$, the sink map 
    has the form $\rad P\to P$.
  \item[2b)] For each indecomposable injective object 
    $I$ in $\fac_\Lambda{\mathcal P}$, the source map has the form 
    $I\to I/\soc I$.
    
    \smallskip
  \item[3a)] There are sufficiently many 
    projective objects in 
    $\sub_\Lambda{\mathcal P}$.
  \item[3b)] There are sufficiently many injective 
    objects in 
    $\fac_\Lambda{\mathcal P}$. \qed
  \end{enumerate}
\end{lemma}

\medskip
In order to obtain information about the injective objects in
$\sub_\Lambda{\mathcal P}$ and the projective objects in 
$\fac_\Lambda{\mathcal P}$ we use that the two categories are in fact
equivalent.

\begin{lemma} \label{subfac}
  There is an equivalence
  of categories between $\sub_\Lambda {\mathcal P}$
  and $\fac_\Lambda{\mathcal P}$. 
\end{lemma}

\begin{proof}
  The equivalence is given by functors
  $$E: \sub_\Lambda {\mathcal P}\to \fac_\Lambda{\mathcal P}\quad\text{and}\quad
  E':\fac_\Lambda {\mathcal P}\to \sub_\Lambda{\mathcal P}$$
  where $E$ maps a system $(M_*,(M_i)_{i\in{\mathcal P}})$ 
  of submodules of $M_*$ to the system 
  $(M_0,(M_0\to\!\!\!\!\!\to M_0/M_i)_{i\in{\mathcal P}})$ of the cokernel maps where
  $M_0=M_*$. Conversely, $E'$ sends a system 
  $(M_0,(M_0\to\!\!\!\!\!\to M_i)_{i\in{\mathcal P}})$ to the object
  $(M_*,\ker(M_0\to M_i)_{i\in{\mathcal P}})$ given by the kernel maps where $M_*=M_0$.
\end{proof}

\medskip Using this equivalence, we obtain the structure of the injective objects
in $\sub_\Lambda\mathcal P$.
We omit the corresponding assertions about the projective objects in 
$\fac_\Lambda\mathcal P$.

\begin{corollary}\label{injectives}
  \smallskip \begin{enumerate}\item[1.] Each injective object in $\sub_\Lambda{\mathcal P}$
    is isomorphic to $E'(I)$ for some injective 
    $\Lambda{\mathcal P}^0$-module $I$.
  \item[2.] For $I$ an indecomposable injective $\Lambda{\mathcal P}^0$-module,
    the morphism $E'(I)\to E'(I/\soc I)$ 
    is a source map in the category $\sub_\Lambda{\mathcal P}^*$.
  \item[3.] There are sufficiently many injective objects in $\sub_\Lambda{\mathcal P}$.
  \end{enumerate}
\end{corollary}

\medskip Here is an even more explicit description of the indecomposable projective and the
indecomposable injective representations in $\sub_\Lambda{\mathcal P}$, and of 
their respective sink and source maps.  

\begin{notation}
  For $M$ a $\Lambda$-module and $S$ a convex subset
  of ${\mathcal P}^*$, define the $\mathcal P^*\Lambda$-module $M^S$ by
  $$M^S_i=\left\{\begin{array}{ll}M&{\rm if}\quad i\in S\\ 0&{\rm otherwise,}\end{array}
  \right.$$
  and let all multiplication maps $M_i^S\to M_j^S, m\mapsto mq_{ij}$ be zero unless
  both $i,j\in S$ and the map can be taken to be $1_M$.
  In particular if $\rho\subset \mathcal P^*\times\mathcal P^*$ is a binary relation and 
  $x\in\mathcal P^*$ then $M^{\rho x}=M^S$ where $S=\{i:(i,x)\in\rho\}$.  Thus, for example
  $M^{\geq x}$ given by
  $$M^{\geq x}_i=\left\{\begin{array}{ll}M&{\rm if}\quad i\geq x\\ 0&{\rm otherwise}\end{array}
  \right.$$
  is a submodule representation of $\mathcal P$.
\end{notation}

\begin{proposition} \label{proj-inj}
  
  \begin{enumerate}\item[1.] If $P$ is an indecomposable
    projective $\Lambda$-mo\-dule and $x\in{\mathcal P}^*$
    then $P^{\geq x}$ is an indecomposable projective
    object in $\sub_\Lambda{\mathcal P}$.  
    The corresponding sink map is the inclusion
    $(\rad P)^{\geq x}+ P^{>x}\subset P^{\geq x}$.
    Each indecomposable projective object 
    in $\sub_\Lambda{\mathcal P}$ has this form.
    
    \smallskip\item[2a)] If $I$ is an indecomposable injective
    $\Lambda$-module and $x\in{\mathcal P}$, then $I^{\not\leq x}$
    is an indecomposable injective object in $\sub_\Lambda{\mathcal P}$.
    The corresponding source map is the inclusion
    $I^{\not\leq x}\subset (\soc I)^{\not< x}+I^{\not\leq x}$.
  \item[2b)] For $I$ an indecomposable injective $\Lambda$-module,
    the representation $I^{\leq*}$ is indecomposable injective
    in $\sub_\Lambda{\mathcal P}$ and has as source map the canonical map
    $I^{\leq*}\to (I/\soc I)^{\leq*}$.
  \item[2c)] Each indecomposable injective object in $\sub_\Lambda{\mathcal P}$
    is isomorphic to one listed above. \qed
  \end{enumerate}
\end{proposition}

\begin{observation} \label{observation}
  In each case the starting term $A$ of a sink map
  for an indecomposable projective object in $\sub_\Lambda\mathcal P$,
  and the end term $C$ of the source map for an indecomposable injective
  object in $\sub_\Lambda\mathcal P$, is indecomposable or zero.
  Either module $A$ or $C$ is determined uniquely by its type. 
\end{observation}

\begin{example}[Chains of Subgroups, I] \label{example-chains}
  We compute the projective and the injective indecomposables in for the category
  $\mathcal C_3(\Lambda)=\sub_\Lambda\mathcal P$ 
  of submodule representations of the linear poset 
  $\mathcal P=\vcenter{\beginpicture 
    \setcoordinatesystem units <.3cm,.3cm>
    \multiput{$\scriptscriptstyle\bullet$} at 0 0  0 1  0 2 /
    \put{$\scriptscriptstyle 3$} at .5 2
    \put{$\scriptscriptstyle 2$} at .5 1
    \put{$\scriptscriptstyle 1$} at .5 0
    \plot 0 0  0 2 /
    \endpicture}$
  where $\Lambda$ is any commutative uniserial ring.
  Thus we are dealing with chains $A=(A_1\subset A_2\subset A_3\subset A_*)$ of 
  submodules of $\Lambda$-modules. The big module $A_*$ is a direct sum of cyclic
  $\Lambda$-modules and hence is given by the partition which lists the lengths of the
  summands.  We picture $A$ by rotating the Young diagram
  of the partition for $A_*$ by $90^{\rm o}$ and by using the symbols 
  $\vcenter{\beginpicture 
    \setcoordinatesystem units <.7cm,.7cm>
    \put{$\smallsq5$} at 0 0
    \multiput{} at 0 -.25  .5 .25 /
    \put{$\mouthm$} at .25 0
    \endpicture}$,
  $\vcenter{\beginpicture 
    \setcoordinatesystem units <.7cm,.7cm>
    \put{$\smallsq5$} at 0 0
    \multiput{} at 0 -.25  .5 .25 /
    \put{$\nosem$} at .25 0
    \endpicture}$, and 
  $\vcenter{\beginpicture 
    \setcoordinatesystem units <.7cm,.7cm>
    \put{$\smallsq5$} at 0 0
    \multiput{} at 0 -.25  .5 .25 /
    \put{$\eyesm$} at .25 0
    \endpicture}$
  to indicate in which of the radical layers in $A_*$ the generators of $A_1$, $A_2$, 
  and $A_3$, respectively, can be found.
  For example if $\Lambda=\mathbb Z/p^2$,
  $A_3=A_*=\Lambda$ and $A_1=A_2=p\mathbb Z/p^2$, then we picture
  $A$ as follows.
  $$A:\quad\vcenter{\beginpicture 
    \setcoordinatesystem units <.7cm,.7cm>
    \multiput{$\smallsq5$} at 0 0  0 .5 /
    \multiput{} at 0 -.5  .5 1 /
    \put{$\eyesm$} at .25 .5
    \put{$\mouthm$} at .25 0 
    \put{$\nosem$} at .25 0
    \endpicture}$$
  For this particular poset, but independent of the length $n$ of $\Lambda$, 
  the projective representation $\Lambda^{\geq x}$ corresponding to a point 
  $x\in\mathcal P^*$ and the injective representation $\Lambda^{\not< x}$
  corresponding to this point coincide. 
  We assume that in addition the length of $\Lambda$ is $n=2$; then the 
  four indecomposable projective-injective objects can be pictured as follows.
  $$\Lambda^{\geq*}:\quad\vcenter{\beginpicture 
    \setcoordinatesystem units <.7cm,.7cm>
    \multiput{$\smallsq5$} at 0 0  0 .5 /
    \multiput{} at 0 -.5  .5 1 /
    \endpicture}\qquad
  \Lambda^{\geq3}:\quad\vcenter{\beginpicture 
    \setcoordinatesystem units <.7cm,.7cm>
    \multiput{$\smallsq5$} at 0 0  0 .5 /
    \multiput{} at 0 -.5  .5 1 /
    \put{$\eyesm$} at .25 .5
    \endpicture}\qquad
  \Lambda^{\geq2}:\quad\vcenter{\beginpicture 
    \setcoordinatesystem units <.7cm,.7cm>
    \multiput{$\smallsq5$} at 0 0  0 .5 /
    \multiput{} at 0 -.5  .5 1 /
    \put{$\eyesm$} at .25 .5
    \put{$\nosem$} at .25 .5
    \endpicture}\qquad
  \Lambda^{\geq1}:\quad\vcenter{\beginpicture 
    \setcoordinatesystem units <.7cm,.7cm>
    \multiput{$\smallsq5$} at 0 0  0 .5 /
    \multiput{} at 0 -.5  .5 1 /
    \put{$\eyesm$} at .25 .5
    \put{$\mouthm$} at .25 .5 
    \put{$\nosem$} at .25 .5
    \endpicture}$$

In this case for $x\in\mathcal P$,
  the radical of the projective presentation $\Lambda^{\geq x}$ coincides with the
  last term of the source map for the injective presentation $\Lambda^{\not\leq x}$.  
  We obtain the
  following part of the Auslander-Reiten quiver for $\sub_\Lambda\mathcal P$.
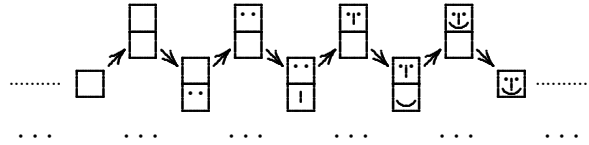
  $$\vcenter{\beginpicture
    \setcoordinatesystem units <.7cm,.7cm> 
    \put{} at  -2 -1
    \put{} at  10 1.5
    \multiput{$\oneo$} at  0 0  8 0 /
    \multiput{$\oo$} at 1 1  2 0  3 1  4 0  5 1  6 0  7 1  /
    \multiput{\eyesm} at  8 0 /
    \multiput{\eyesl} at 2 0 /
    \multiput{\eyesh} at  3 1  4 0  5 1  6 0  7 1 /
    \multiput{\nosem} at 8 0 /
    \multiput{\nosel} at 4 0 /
    \multiput{\noseh} at 5 1  6 0  7 1 /
    \multiput{\mouthm} at 8 0 /
    \multiput{\mouthl} at 6 0 /
    \multiput{\mouthh} at 7 1 /
    \arr{0.35 0.35} {0.65 0.65}
    \arr{2.35 0.35} {2.65 0.65}
    \arr{4.35 0.35} {4.65 0.65}
    \arr{6.35 0.35} {6.65 0.65}
    \arr{1.35 0.65} {1.65 0.35}
    \arr{3.35 0.65} {3.65 0.35}
    \arr{5.35 0.65} {5.65 0.35}
    \arr{7.35 0.65} {7.65 0.35}
    
    \setdots<2pt>
    \plot -1.5 0  -.5 0 /
    \plot 8.5 0  9.5 0 /
    \multiput{$\cdots$} at -1 -1  1 -1  3 -1  5 -1  7 -1  9 -1 /
    \endpicture}
  $$
\end{example}

\section{Auslander-Reiten sequences 
  which are not split exact in each
  component}
\label{section-nonsplit}

Suppose that $\mathcal E:0\to A\to B\to C\to 0$ is an Auslander-Reiten sequence in the submodule
category $\sub_\Lambda\mathcal P$. For $x\in\mathcal P^*$ one can consider the 
short exact sequence $\mathcal E_x:0\to A_x\to B_x\to C_x\to 0$ consisting of the $x$-components.
It turns out that $\mathcal E_x$ is either split exact or almost split. We show that given $\mathcal E$,
then at most one of the sequences $\mathcal E_x$ is not split exact.
Moreover, almost all Auslander-Reiten sequences are split exact in every component. 

\smallskip We first determine the exceptions.
Let ${\mathcal T}:0\to U\sto rV\sto sW\to 0$ 
be an Auslander-Reiten sequence in the category
of modules over the base ring $\Lambda$, and let $x\in{\mathcal P}^*$.
Define the sequence 
$${\mathcal E}={\mathcal E}({\mathcal T},x): \quad 0\longrightarrow A\lto fB\lto gC
\longrightarrow 0$$
as follows.
Let $i:U\to E$ be the inclusion of $U$ in its injective envelope,
and let $j:V\to E$ be a lifting of $i$ over the monomorphism $r:U\to V$,
so $jr=i$ holds. Then there are modules
$A,B,C\in\sub_\Lambda{\mathcal P}$ given by the following $\Lambda$-modules where
$y\in\mathcal P^*$.
$$ A_y=\left\{\begin{matrix} U, & y\leq x \cr E, & y\not\leq x\end{matrix}\right.,
\quad B_y=\left\{\begin{matrix} V, & y=x \cr A_y\oplus C_y, & y\neq x
  \end{matrix}\right., \quad
C_y=\left\{\begin{matrix} W, & y\geq x \cr 0, & y\not\geq x \end{matrix}\right.$$
It is clear which inclusion maps make up the modules $A$ and $C$;
for $B$, the subspaces $E$, $U$, $V$ are embedded in $E\oplus W$
via the maps ${1\choose 0}:E\to E\oplus W$, $i:U\to E$, $r:U\to V$,
and ${j\choose s}:V\to E\oplus W$.   

\smallskip The maps $f:A\to B$ and $g:B\to C$ are given as follows.
For $y\neq x$, $f_y$ and $g_y$ are the canonical inclusions and projections,
respectively, while $f_x=r$ and $g_x=s$.

\smallskip 
Note that there are at most 4 possiblities for the components $\mathcal E_y$ of $\mathcal E$,
depending on whether $y=x$, $y<x$, $y>x$ or $y$ and $x$ are unrelated.
The following example is such that all 4 possibilities are realized.

\begin{example}
  If the poset is 
  ${\mathcal P}^*=\vcenter{\beginpicture 
    \setcoordinatesystem units <.4cm,.4cm>
    \multiput{$\scriptscriptstyle\bullet$} at 0 0  1 -1  2 0 /
    \put{$\scriptstyle*$} at 1 1
    \put{$\scriptscriptstyle 2$} at -.3 0
    \put{$\scriptscriptstyle 3$} at 2.3 0
    \put{$\scriptscriptstyle 1$} at 1.3 -1
    \plot .3 .3  .7 .7 /
    \plot 1.3 .7  1.7 .3 /
    \plot .3 -.3  .7 -.7 /
    \plot 1.3 -.7  1.7 -.3 /
    \endpicture}$, and $x=2$, then the Auslander-Reiten sequence 
  $\mathcal T:0\to U\to^r V\to^sW\to 0$ in $\mod \Lambda$
  gives rise to the Auslander-Reiten sequence $\mathcal E$ in $\sub_\Lambda\mathcal P$:
  $${\mathcal E}:\quad 
  0 
  \lto{} 
  \quad 
  \prep EUEU{i}{1}{1}{i}
  \quad \lto{\pmap{({\scriptscriptstyle 1\atop \scriptscriptstyle0})}r11} \quad
  \prep {E\oplus W}VEU{{j\choose s}}{{1\choose0}}{r}{i} 
  \quad \lto{\pmap{(0\;1)}s00} \quad
  \prep WW001{}{}{}
  \lto{} 0$$ 
\end{example}

\begin{theorem} \label{almost-split}
  \begin{enumerate}\item[1.] Let ${\mathcal T}$ be an almost split
    sequence in $\mod \Lambda$ and $x\in {\mathcal P}^*$. 
    Then the sequence 
    ${\mathcal E}={\mathcal E}({\mathcal T},x)$ is an Auslander-Reiten sequence
    in $\sub_\Lambda{\mathcal P}$ which is split exact in every component except in
    the $x$-component.
    
  \item[2.] Conversely if ${\mathcal E}:0\to A\to B\to C\to 0$ is an Auslander-Reiten sequence
    in $\sub_\Lambda{\mathcal P}$, and $x\in{\mathcal P}^*$ is such
    that the sequence in $\mod\Lambda$ of the $x$-components
    ${\mathcal T}:0\to A_x\to B_x\to C_x\to 0$ is not split
    exact, then $\mathcal T$ is an Auslander-Reiten sequence
    in $\mod\Lambda$ and the sequences ${\mathcal E}$ and 
    ${\mathcal E}({\mathcal T},x)$ are equivalent.
  \end{enumerate}
\end{theorem}

\begin{proof}
  Clearly, the sequence ${\mathcal E}:0\to A\to B\to C\to 0$ 
  in the first assertion is not split exact
  and the modules $A$ and $C$ are indecomposable.  We show that ${\mathcal E}$
  is right almost split.
  
  \smallskip Let $T\in\sub_\Lambda{\mathcal P}$ and let $t:T\to C$ be a morphism
  which is not a split epimorphism.  Then $t_x:T_x\to C_x$ is not a split
  epimorphism in $\mod\Lambda$ and hence factors over 
  $g_x:V\to W$: There is $t_x':T_x\to B_x$ such that $g_xt_x'=t_x$, as illustrated below in
  the situation of the example.
  $$
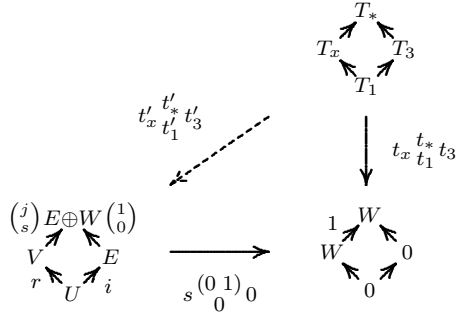
\beginpicture\linethickness1mm
  \setcoordinatesystem units <1.3cm,.9cm>
  \put{$\prep{E\oplus W}VEU{{j\choose s}}{{1\choose0}}{r}{i}$} at 0 0
  \put{$\prep WW00{1}{}{}{}$} at 3 0
  \put{$\prep {T_*}{T_x}{T_3}{T_1}{}{}{}{}$} at 3 3
  \arr{1 0}{2 0}
  \put{$\pmap{(0\;1)}s00$} at 1.5 -.6
  \arr{3 2}{3 1}
  \put{$\pmap{t_*}{t_x}{t_3}{t_1}$} at 3.6 1.5
  \setdashes<2pt>
  \arr{2 2}{1 1}
  \put{$\pmap{t_*'}{t_x'}{t_3'}{t_1'}$} at 1 2
  \endpicture$$
  If $x\neq *$, define $t_*':T_*\to E\oplus W$ as follows.  Since $E$ is injective and 
  $T_x\subseteq T_*$, the map $jt_x':T_x\to E$ extends to a map
  $t_*^+:T_*\to E$, hence  $t_*^+|_{T_x}=jt_x'$ holds.
  Put $t_*'={t_*^+ \choose t_*}$.
  As a consequence, for $y\geq x$, the map $t_y':T_y\to E\oplus W$
  is defined uniquely as the restriction of $t_*'$ to $T_y$. 
  
  Next for $y\in{\mathcal P}$ which is not in relation with $x$,
  define the map $t_y':T_y\to E$ as follows.
  Since $t_*|_{T_y}=0$, the image $t_*'(T_y)$ is contained in $E\oplus0$
  and hence there is a unique map $t_y':T_y\to E$ such that 
  ${1\choose0}t_y'=t_*'|_{T_y}: T_y\to E\oplus W$.
  
  Finally consider $y\in{\mathcal P}$ such that $y<x$.
  Since $t_x|_{T_y}=0$, also $st_x'|_{T_y}=0$ and hence the restriction
  of $t_x'$to $T_y$ factors over $r:U\to V$: There exists a unique
  $t_y': T_y\to U$ such that $t_x'|_{T_y}=rt_y'$.
  
  In conclusion, $t'=(t'_y)_{y\in{\mathcal P}}: T\to B$ is a morphism 
  which satisfies $gt'=t$.

  \medskip
  It remains to show the second assertion.  Note that all the sequences constructed
  have the property that the last term has the form $W^{\geq x}$ where $W$ is an
  indecomposable $\Lambda$-module and $x\in{\mathcal P}^*$.  We show that 
  any Auslander-Reiten sequence ${\mathcal E}:0\to A\sto fB\sto gC\to 0$
  for which the last term is not of type $W^{\geq x}$ for a $\Lambda$-module $W$,
  is split exact in each component.  
  For $y\in{\mathcal P}^*$ consider the sequence of the $y$-components 
  $$0\to A_y\sto{f_y}B_y\sto{g_y}C_y\to 0$$ 
  in $\mathcal E$. 
  Define the $\mathcal P$-representation $T=(C_y)^{\geq y}$ and
  consider the map $t:T\to C$ given by the identity map on the $y$-component.
  By assumption, $C$ has not the form $W^{\geq x}$ so
  $t$ is not an epimorphism and hence not a split epimorphism.  Since $t$  factors over $g$,
  the map $g_y: B_y\to C_y$ must be a split epimorphism in $\mod \Lambda$. 
  This finishes the proof. 
\end{proof}

\begin{example}[Chains of subgroups, II] \label{chains-two}
  We use the notation from Example \ref{example-chains} where 
  $\mathcal P=\vcenter{\beginpicture 
    \setcoordinatesystem units <.3cm,.3cm>
    \multiput{$\scriptscriptstyle\bullet$} at 0 0  0 1  0 2 /
    \put{$\scriptscriptstyle 3$} at .5 2
    \put{$\scriptscriptstyle 2$} at .5 1
    \put{$\scriptscriptstyle 1$} at .5 0
    \plot 0 0  0 2 /
    \endpicture}$ and $\Lambda$ is a commutative uniserial ring of length 2.
  Denote the only Auslander-Reiten sequence in $\mod\Lambda$ by
  $\mathcal T:0\to k\sto\iota \Lambda\sto\pi k\to 0$ where $k=\Lambda/\rad\Lambda$.
  For $x=*$ we obtain the Auslander-Reiten sequence in $\sub_\Lambda\mathcal P$:
  $$\mathcal E(\mathcal T,*): 
  \quad 0\;\to\; ({\scriptstyle k=k=k=k})\; \to\; ({\scriptstyle k=k=k\sto\iota \Lambda})
  \;\to\; ({\scriptstyle 0=0=0\subset k})\;\to 0$$
  For $x\in\mathcal P$, the sequence $\mathcal E(\mathcal T,x)$ which we picture
  here for $x=2$
  \begin{eqnarray*}
    \mathcal E(\mathcal T,2): \\
    0\;\to & & ({\scriptstyle k=k\sto\iota\Lambda=\Lambda})
    \lto{1,\iota,{1\choose0},{1\choose0}}
    ({\scriptstyle k\sto\iota\Lambda\lto{{1\choose\pi}}(\Lambda\oplus k)=(\Lambda\oplus k)}) \lto{}\\
    & & \lto{0,\pi,{\scriptscriptstyle (0\,1)},{\scriptscriptstyle (0\,1)}}
    ({\scriptstyle 0\subset k=k=k})\;\to\; 0
  \end{eqnarray*}
  has decomposable middle term:
  $$ ({\scriptstyle k\sto\iota\Lambda\lto{{1\choose\pi}}(\Lambda\oplus k)=(\Lambda\oplus k)})
  \quad  \cong \quad 
  ({\scriptstyle k\sto\iota\Lambda=\Lambda=\Lambda})\quad\oplus
  \quad({\scriptstyle 0=0\subset k=k}).$$
  Note that the irreducible map onto the second summand of the middle term 
  of $\mathcal E(\mathcal T,x)$ coincides with the irreducible map 
  starting at the first summand of the middle term in the sequence
  $\mathcal E(\mathcal T,x+1)$ (in case $x=3$, take the right almost split map in 
  $\mathcal E(\mathcal T,*)$). Thus the Auslander-Reiten sequences $\mathcal E(\mathcal T,x)$
  for $x=*,3,2,1$ line up to form the following part of the
  Auslander-Reiten quiver for $\sub_\Lambda\mathcal P$:
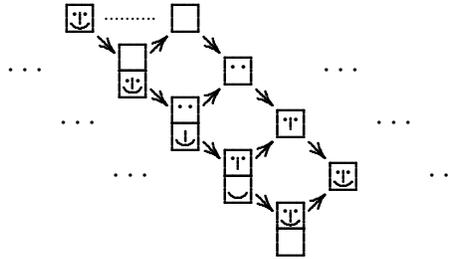
  $$\vcenter{\beginpicture
    \setcoordinatesystem units <.7cm,.7cm> 
    \put{} at  -1 5
    \put{} at  7 -1
    \multiput{$\oneo$} at  0 4  2 4  3 3  4 2  5 1 /
    \multiput{$\oo$} at 1 3  2 2  3 1  4 0 /
    \multiput{\eyesm} at  0 4  3 3  4 2  5 1 /
    \multiput{\eyesl} at 1 3 /
    \multiput{\eyesh} at 2 2  3 1  4 0 /
    \multiput{\nosem} at 0 4  4 2  5 1 /
    \multiput{\nosel} at 1 3  2 2  /
    \multiput{\noseh} at 3 1  4 0 /
    \multiput{\mouthm} at 0 4  5 1 /
    \multiput{\mouthl} at 1 3  2 2  3 1  /
    \multiput{\mouthh} at 4 0 /
    \arr{1.35 3.35} {1.65 3.65}
    \arr{2.35 2.35} {2.65 2.65}
    \arr{3.35 1.35} {3.65 1.65}
    \arr{4.35 0.35} {4.65 0.65}
    \arr{0.35 3.65} {0.65 3.35}
    \arr{1.35 2.65} {1.65 2.35}
    \arr{2.35 1.65} {2.65 1.35}
    \arr{3.35 0.65} {3.65 0.35}
    \arr{2.35 3.65} {2.65 3.35}
    \arr{3.35 2.65} {3.65 2.35}
    \arr{4.35 1.65} {4.65 1.35}

    \setdots<2pt>
    \plot .5 4  1.5 4 /
    \multiput{$\cdots$} at  -1 3  0 2  1 1   7 1 6 2  5 3  /  
    \endpicture}
  $$

\end{example}

\medskip
There is also the following corresponding result for the category
$\fac_\Lambda{\mathcal P}$.  Again, let $0\to U\sto rV\sto sW\to 0$
be an Auslander-Reiten sequence in the category $\mod \Lambda$ and
$x\in{\mathcal P}^0$. Let $p:P\to W$ be a projective cover and 
$q:P\to V$ a lifting of $p$ over the epimorphism $s:V\to W$.
Define the representations $A$, $B$, $C$ as follows 
$$ A_y=\left\{\begin{matrix} U, & y\leq x \cr 0, & y\not\leq x\end{matrix}\right.,
\quad
B_y=\left\{\begin{matrix} V, & y=x\cr A_y\oplus C_y, & y\neq x
  \end{matrix}\right.,
\quad C_y=\left\{\begin{matrix} W, & y\geq x \cr P, & y\not\geq x \end{matrix}\right.,$$
where $y\in\mathcal P^0$.
\begin{example}
  We illustrate the linear maps which make up $A$, $B$, $C$, $f$, and $g$
  in the example where
  ${\mathcal P}^0=\vcenter{\beginpicture 
    \setcoordinatesystem units <.4cm,.4cm>
    \multiput{$\scriptscriptstyle\bullet$} at 0 0  1 1  2 0 /
    \put{$\scriptstyle0$} at 1 -1
    \put{$\scriptscriptstyle 1$} at -.3 0
    \put{$\scriptscriptstyle 2$} at 2.3 0
    \put{$\scriptscriptstyle 3$} at 1.3 1
    \plot .3 .3  .7 .7 /
    \plot 1.3 .7  1.7 .3 /
    \plot .3 -.3  .7 -.7 /
    \plot 1.3 -.7  1.7 -.3 /
    \endpicture}$, and $x=1$.
  The Auslander-Reiten sequence in $\fac_\Lambda{\mathcal P}$ is as follows.
  $${\mathcal E}:\quad 
  0 
  \lto{}
  \;
  \prep 0U0U{}{}1{} 
  \quad \lto{\pmap0r0{({\scriptscriptstyle 1\atop \scriptscriptstyle0})}} \quad
  \prep WVP{U\oplus P}sp{(r\;q)}{(0\;1)}
  \quad \lto{\pmap1s1{(0\;1)}} \quad
  \prep WWPP1pp1
  \lto{} 0$$ 
\end{example}

\begin{theorem} \label{almost-split-factors}
  \begin{enumerate}\item[1.] Let ${\mathcal T}$ be an almost split
    sequence in $\mod \Lambda$ and $x\in {\mathcal P}^0$. 
    Then the sequence 
    ${\mathcal E}={\mathcal E}({\mathcal T},x)$ is an Auslander-Reiten sequence
    in $\fac_\Lambda{\mathcal P}$ which is split exact in every component 
    except in the $x$-component.
    
  \item[2.] Conversely if ${\mathcal E}:0\to A\to B\to C\to 0$ is an Auslander-Reiten sequence
    in $\fac_\Lambda{\mathcal P}$, and $x\in{\mathcal P}^0$ is such
    that the sequence in $\mod\Lambda$ of the $x$-components
    ${\mathcal T}:0\to A_x\to B_x\to C_x\to 0$ is not split
    exact, then $\mathcal T$ is an Auslander-Reiten sequence in $\mod\Lambda$ and the
    sequences ${\mathcal E}$ and ${\mathcal E}({\mathcal T},x)$ are equivalent.\qed
  \end{enumerate}
\end{theorem}

\begin{example}[Bounded submodules of modules, I] \label{bounded-I}
  In this example we address the remark by M.C.R.~Butler.
  For $\mathcal P=\bullet_1$ the one point poset, write $\mathcal S(\Lambda)=
  \sub_\Lambda\mathcal P$ and let for $m\leq n$ the full subcategory
  $\mathcal S_m(\Lambda)$ of $\mathcal S(\Lambda)$ consist of all pairs 
  $(A_1\subset A_*)$ where $\rad^m\!A_1=0$ holds.  Let $t$ be a radical
  generator
  for $\Lambda$.

  \smallskip
  The category $\mathcal S_m(\Lambda)$ is an exact Krull-Remak-Schmidt category 
  with Auslander-Reiten sequences \cite{as}. Namely, 
  each object $(A_1\subset A_*)$ in $\mathcal S(\Lambda)$
  has a minimal right and a minimal left approximation in $\mathcal S_m(\Lambda)$,
  given as follows.
  $$ \big(\soc^m\!A_1\subset A_*\big)\to \big(A_1\subset A_*\big);$$
  $$\big(A_1\subset A_*\big)\to \big(A_1/\rad^m\!A_1\subset A_*/\rad^m\!A_1\big)$$

  \smallskip We determine each Auslander-Reiten sequence 
  for which one of the short exact sequences given by either (1) the submodules, (2) the 
  total spaces, or (3) the factor modules is not split exact.

  \begin{enumerate} 
    \item The sequence of submodules is not split exact.  
      Let $\bar\Lambda$ be the factor ring $\Lambda/\rad^m\Lambda$ and suppose that 
      $0\to A_1\to B_1\to C_1\to 0$ is an Auslander-Reiten sequence in $\mod\bar\Lambda$.
      Let $u:A_1\to E$ be an injective envelope in $\mod\Lambda$ and 
      choose an extension $v:B_1\to E$ of $u$.  According to Theorem \ref{almost-split},
      the sequence
      $$0\to \big(A_1\to E\big)\to \big(B_1\to E\oplus C_1)\to \big(C_1=C_1\big)\to 0$$
      is an Auslander-Reiten sequence in $\mathcal S(\Lambda)$;
      since all objects are in $\mathcal S_m(\Lambda)$, the sequence is an 
      Auslander-Reiten sequence in $\mathcal S_m(\Lambda)$, too.
      For example if $n=6$ and $m\geq3$ then
      the Auslander-Reiten sequence in $\mathcal S_m(\Lambda)$ ending at 
      $C$ where $C_1=C_*=\Lambda/\rad^2\!\Lambda$ is as follows.
      $$0 \; \longrightarrow \;
      \vcenter{
        \beginpicture 
        \setcoordinatesystem units <0.25cm,0.25cm>
        \multiput{} at 0 -2  1 6 /
        \multiput{\sq} at 0 5  0 4  0 3  0 2  0 1  0 0 /
        \put{$\scriptstyle \bullet$} at 0.5 1 
        \endpicture}
      \;\longrightarrow \;
      \vcenter{
        \beginpicture 
        \setcoordinatesystem units <0.25cm,0.25cm>
        \multiput{} at 0 -2 1 6 /
        \multiput{\sq} at 0 5  0 4  0 3  0 2  0 1  0 0  1 2  1 1 /
        \put{$\scriptstyle \bullet$} at 0.5 2
        \put{$\scriptstyle \bullet$} at 1.5 2 
        \put{$\scriptstyle \bullet$} at 1.5 1 
        \plot 0.5 2  1.5 2 /
        \endpicture}
      \;\longrightarrow \;
      \vcenter{\beginpicture 
        \setcoordinatesystem units <0.25cm,0.25cm>
        \multiput{} at 0 -2  1 6 /
        \multiput{\sq} at 0 2  0 1 /
        \put{$\scriptstyle \bullet$} at 0.5 2 
        \endpicture}
      \;\longrightarrow \; 0 $$
      The objects have type $t_A=(2;6)$, $t_B=(31;62)$, and $t_C=(2;2)$.
      As a consequence, the middle term of this sequence is an indecomposable object.  
      It is given by the inclusion $(B_1\subset B_*)$
      where $B_*$ is generated by two elements $y_1$ and $y_2$, say, bounded by $t^6$ and
      $t^2$, respectively.  The submodule $B_1$ also has two generators, for example
      $x_1=t^3y_1+y_2$ (indicated in the diagram by the connected dots) 
      and $x_2=ty_2$ (given by the isolated point).
    \item The sequence of the total spaces is not split exact.  
      Let $0\to A_*\to B_*\to C_*\to 0$ be an Auslander-Reiten sequence in $\mod\Lambda$.
      According to Theorem \ref{almost-split}, 
      $$0\to \big(A_*=A_*\big)\to \big(A_*\subset B_*\big)\to \big(0\subset C_*\big)
      \to 0$$
      is an Auslander-Reiten sequence in $\mathcal S(\Lambda)$ and, using right
      approximations, 
      $$\qquad 0\to \big(\soc^m\!A_*\subset A_*\big)\to \big(\soc^m\!A_*\subset B_*\big)
      \to \big(0\subset C_*\big)\to 0$$
      is an Auslander-Reiten sequence in $\mathcal S_m(\Lambda)$. For example 
      if $n\geq 5$,
      and $m=3$, then the Auslander-Reiten sequence in
      $\mathcal S_3(\Lambda)$ starting at $A$ where $A_*=\Lambda/\rad^4\!\Lambda$
      and $A_1=\soc^3\!A_*$ is as follows.
      $$0 \; \longrightarrow \;
      \vcenter{
        \beginpicture 
        \setcoordinatesystem units <0.25cm,0.25cm>
        \multiput{} at 0 -2  1 5 /
        \multiput{\sq} at 0 3  0 2  0 1  0 0 /
        \put{$\scriptstyle \bullet$} at 0.5 2 
        \endpicture}
      \;\longrightarrow \;
      \vcenter{
        \beginpicture 
        \setcoordinatesystem units <0.25cm,0.25cm>
        \multiput{} at 0 -2 1 5 /
        \multiput{\sq} at 1 4  0 3  0 2  0 1  1 0  1 3  1 2  1 1 /
        \put{$\scriptstyle \bullet$} at 0.5 2
        \put{$\scriptstyle \bullet$} at 1.5 2 
        \plot 0.5 2  1.5 2 /
        \endpicture}
      \;\longrightarrow \;
      \vcenter{\beginpicture 
        \setcoordinatesystem units <0.25cm,0.25cm>
        \multiput{} at 0 -2  1 5 /
        \multiput{\sq} at 0 4  0 3  0 2  0 1 /
        \endpicture}
      \;\longrightarrow \; 0 $$
      The modules have type $t_A=(3;4)$, $t_B=(3;53)$, and $t_C=(-;4)$; the middle term
      $B$ is an indecomposable object. 
    \item The sequence of the factor modules is not split exact.  

      In this section we consider sequences $0\to A\to B\to C\to 0$ in $\mathcal S(\Lambda)$
      which are split exact in every component, but for which the (short exact) sequence of factors
      $0\to A_*/A_1\to B_*/B_1\to C_*/C_1\to 0$ is almost split. 
      Let 
      $0\to \bar A\sto f\bar B\sto g\bar C\to 0$ be an Auslander-Reiten sequence in 
      $\mod\Lambda$, $u:P\to \bar C$ a projective cover with kernel $v:\bar K\to P$
      and $w:P\to \bar B$ a lifting of $u$.  Let $z:K\to \bar A\oplus P$ be the kernel
      of $(f,w)$.  By Theorem \ref{almost-split-factors}, we have the 
      Auslander-Reiten sequence in $\fac_\Lambda\mathcal P$.
      $$0\to \big( \bar A=\bar A\big) \to \big( \bar A\oplus P\to \bar B\big) \to
      \big(P\to \bar C\big)\to 0 $$
      Using the equivalence between $\sub_\Lambda\mathcal P$ and $\fac_\Lambda\mathcal P$
      given by the kernel functor (Lemma~\ref{subfac}), 
      we obtain the following Auslander-Reiten sequence
      in $\mathcal S(\Lambda)$.
      $$0\to \big( 0\subset \bar A\big) \to \big( K\to \bar A\oplus P\big) \to
      \big(\bar K\to P\big)\to 0$$
      (By the snake lemma, the two kernels $K$ and $\bar K$ are isomorphic.) 
      Using left approximations, we arrive at the Auslander-Reiten sequence in
      $\mathcal S_m(\Lambda)$.
      \begin{eqnarray*} \lefteqn{\hspace{7ex}0\to \big( 0\subset \bar A\big) 
          \to \big( K/\rad^m\!K\to (\bar A\oplus P)/z(\rad^m\!K) \big) \to }  \\
        & & \hspace{25ex}\to (\bar K/\rad^m\bar K\to P/v(\rad^m\!\bar K)\big)\to 0\\
      \end{eqnarray*}
      For example if $n\geq 6$, $m=3$, and $\bar A=\Lambda/\rad^2\!\Lambda$ we obtain
      the following Auslander-Reiten sequence in the category $\mathcal S_3(\Lambda)$.
      If the {\it cotype} of an object $A$ is the type of the factor $A_*/A_1$ then
      the modules have cotypes $\cot_A=(2)$, $\cot_B=(31)$, and $\cot_C=(2)$. In particular,
      $B$ is an indecomposable object.
      $$0 \; \longrightarrow \;
      \vcenter{
        \beginpicture 
        \setcoordinatesystem units <0.25cm,0.25cm>
        \multiput{} at 0 -1  1 5 /
        \multiput{\sq} at 0 3  0 2  /
        \endpicture}
      \;\longrightarrow \;
      \vcenter{
        \beginpicture 
        \setcoordinatesystem units <0.25cm,0.25cm>
        \multiput{} at 0 -1 1 5 /
        \multiput{\sq} at 1 4  0 3  0 2  1 1  1 0  1 3  1 2  /
        \put{$\scriptstyle \bullet$} at 0.5 2
        \put{$\scriptstyle \bullet$} at 1.5 2 
        \plot 0.5 2  1.5 2 /
        \endpicture}
      \;\longrightarrow \;
      \vcenter{\beginpicture 
        \setcoordinatesystem units <0.25cm,0.25cm>
        \multiput{} at 0 -1  1 5 /
        \multiput{\sq} at 0 4  0 3  0 2  0 1  0 0 /
        \put{$\scriptstyle \bullet$} at .5 2
        \endpicture}
      \;\longrightarrow \; 0 $$
  \end{enumerate}
  It turns out that in the case $m=3$, $n=6$,
  the first and last modules in the above Auslander-Reiten sequences 
  form one stable orbit under the translation to which the projective injective 
  module is attached, 
  $$
  
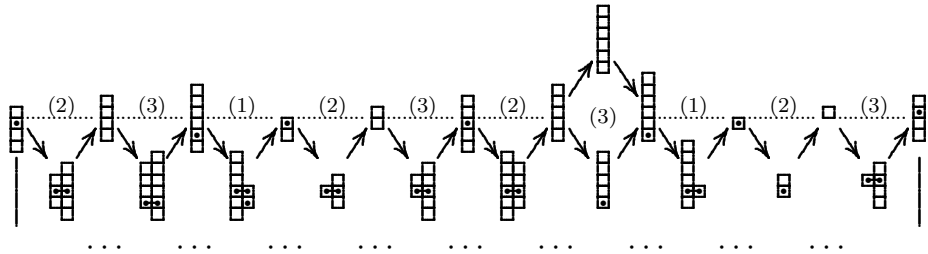
\beginpicture\setcoordinatesystem units <3mm,3mm>
  
  \multiput{$\smallsq5$} at 1.5 16  1.5 16.5  1.5 17
  2 15.5  2 16  2 16.5  2 17  2 17.5 /
  \multiput{\bulletsymb} at 1.75 16.5  2.25 16.5 /
  \plot  1.75 16.5  2.25 16.5 /
  
  \multiput{$\smallsq5$} at 5.5 16  5.5 16.5  5.5 17  5.5 17.5  6 15.5
  6 16  6 16.5  6 17  6 17.5  6 18 /
  \multiput{\bulletsymb} at 5.75 16  6.25 16 /
  \plot 5.75 16  6.25 16 /
  
  \multiput{$\smallsq5$} at 9.5 15.5  9.5 16  9.5 16.5  9.5 17  9.5 17.5
  9.5 18  10 16  10 16.5 /
  \multiput{\bulletsymb} at 9.75 16.5  10.25 16.5  10.25 16 /
  \plot 9.75 16.5  10.25 16.5 /
  
  \multiput{$\smallsq5$} at 13.5 16.5  14 16  14 16.5  14 17 /
  \multiput{\bulletsymb} at 13.75 16.5  14.25 16.5 /
  \plot 13.75 16.5  14.25 16.5 /
  
  \multiput{$\smallsq5$} at 17.5 16.5  17.5 17  18 15.5  18 16  18 16.5
  18 17  18 17.5 /
  \multiput{\bulletsymb} at 17.75 16.5  18.25 16.5 /
  \plot 17.75 16.5  18.25 16.5 /
  
  \multiput{$\smallsq5$} at 21.5 15.5  21.5 16  21.5 16.5  21.5 17  21.5 17.5
  21.5 18  22 16  22 16.5  22 17  22 17.5 /
  \multiput{\bulletsymb} at 21.75 16.5  22.25 16.5 /
  \plot 21.75 16.5  22.25 16.5 /
  
  \multiput{$\smallsq5$} at 25.75 16  25.75 16.5  25.75 17  25.75 17.5
  25.75 18 /
  \multiput{\bulletsymb} at 26 16 /
  
  \multiput{$\smallsq5$} at 29.5 16  29.5 16.5  29.5 17  29.5 17.5  29.5 18
  29.5 18.5  30 16.5 /
  \multiput{\bulletsymb} at 29.75 16.5  30.25 16.5 /
  \plot  29.75 16.5  30.25 16.5 /
  
  \multiput{$\smallsq5$} at 33.75 16.5  33.75 17 /
  \multiput{\bulletsymb} at 34 16.5 /
  
  \multiput{$\smallsq5$} at 37.5 17  38 16  38 16.5  38 17  38 17.5 /
  \multiput{\bulletsymb} at 37.75 17  38.25 17 /
  \plot 37.75 17  38.25 17 /
  
  \multiput{$\smallsq5$} at -.25 18.5  -.25 19  -.25 19.5  -.25 20 /
  \multiput{\bulletsymb} at 0 19.5 /
  
  \multiput{$\smallsq5$} at 3.75 19  3.75 19.5  3.75 20  3.75 20.5 /
  
  \multiput{$\smallsq5$} at 7.75 18.5  7.75 19  7.75 19.5  7.75 20 
  7.75 20.5  7.75 21 /
  \multiput{\bulletsymb} at 8 19 /
  
  \multiput{$\smallsq5$} at 11.75 19  11.75 19.5 /
  \multiput{\bulletsymb} at 12 19.5 /
  
  \multiput{$\smallsq5$} at 15.75 19.5  15.75 20 /
  
  \multiput{$\smallsq5$} at 19.75 18.5  19.75 19  19.75 19.5  19.75 20
  19.75 20.5  /
  \multiput{\bulletsymb} at 20 19.5 /
  
  \multiput{$\smallsq5$} at 23.75 19  23.75 19.5  23.75 20  23.75 20.5
  23.75 21 / 
  
  \multiput{$\smallsq5$} at 27.75 19  27.75 19.5  27.75 20  27.75 20.5
  27.75 21  27.75 21.5 /
  \multiput{\bulletsymb} at 28 19 /
  
  \multiput{$\smallsq5$} at 31.75 19.5 /
  \multiput{\bulletsymb} at 32 19.5 /
  
  \multiput{$\smallsq5$} at 35.75 20 /
  
  \multiput{$\smallsq5$} at 39.75 19  39.75 19.5  39.75 20  39.75 20.5 /
  \multiput{\bulletsymb} at 40 20 /
  \multiput{$\smallsq5$} at 25.75 22  25.75 22.5  25.75 23  25.75 23.5
  25.75 24  25.75 24.5 /
\arr{0.5 19.25} {1.33 18}
\arr{4.5 19.25} {5.33 18}
\arr{8.5 19.25} {9.33 18}
\arr{12.5 19.25} {13.33 18}
\arr{16.5 19.25} {17.33 18}
\arr{20.5 19.25} {21.33 18}
\arr{24.5 19.25} {25.33 18}
\arr{28.5 19.25} {29.33 18}
\arr{32.5 19.25} {33.33 18}
\arr{36.5 19.25} {37.33 18}
\arr{2.67 18} {3.5 19.25}
\arr{6.67 18} {7.5 19.25}
\arr{10.67 18} {11.5 19.25}
\arr{14.67 18} {15.5 19.25}
\arr{18.67 18} {19.5 19.25}
\arr{22.67 18} {23.5 19.25}
\arr{26.67 18} {27.5 19.25}
\arr{30.67 18} {31.5 19.25}
\arr{34.67 18} {35.5 19.25}
\arr{38.67 18} {39.5 19.25}
  \arr{24.5 20.75} {25.5 22.25}
  \arr{26.5 22.25} {27.5 20.75}
  \put{$\scriptstyle(2)$} at 2 20.25
  \put{$\scriptstyle(3)$} at 6 20.25
  \put{$\scriptstyle(1)$} at 10 20.25
  \put{$\scriptstyle(2)$} at 14 20.25
  \put{$\scriptstyle(3)$} at 18 20.25
  \put{$\scriptstyle(2)$} at 22 20.25
  \put{$\scriptstyle(3)$} at 26 19.75
  \put{$\scriptstyle(1)$} at 30 20.25
  \put{$\scriptstyle(2)$} at 34 20.25
  \put{$\scriptstyle(3)$} at 38 20.25
  \setdots<2pt>
  \plot 0.5 19.75  3.5 19.75 /
  \plot 4.5 19.75  7.5 19.75 /
  \plot 8.5 19.75  11.5 19.75 /
  \plot 12.5 19.75  15.5 19.75 /
  \plot 16.5 19.75  19.5 19.75 /
  \plot 20.5 19.75  23.5 19.75 /
  \plot 28.5 19.75  31.5 19.75 /
  \plot 32.5 19.75  35.5 19.75 /
  \plot 36.5 19.75  39.5 19.75 /
  \setsolid
  \plot 0 18  0 15 /
  \plot 40 18  40 15 /
  \multiput{$\cdots$} at 4 14  8 14  12 14  16 14  20 14  24 14  28 14  32 14  36 14 /
  \endpicture$$
  and also an orbit of length three from the projective noninjective
  indecomposable to the injective nonprojective indecomposable, as pictured below.
  $$
  
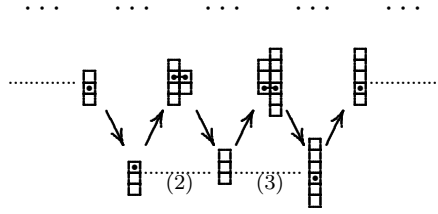
\beginpicture\setcoordinatesystem units <3mm,3mm>
  \multiput{$\smallsq5$} at  25.75 -5  25.75 -4.5  25.75 -4 /
  \multiput{\bulletsymb} at 26 -4 /
  
  \multiput{$\smallsq5$} at 29.75 -4.5  29.75 -4  29.75 -3.5 /
  
  \multiput{$\smallsq5$} at 33.75 -5.5  33.75 -5  33.75 -4.5  33.75 -4
  33.75 -3.5  33.75 -3 /
  \multiput{\bulletsymb} at 34 -4.5 /
  \multiput{$\smallsq5$} at 23.75 -1  23.75 -.5  23.75 0  /
  \multiput{\bulletsymb} at 24 -.5 /
  
  \multiput{$\smallsq5$} at 27.5 -1  27.5 -.5  27.5 0  27.5 .5 28 -.5  28 0 /
  \multiput{\bulletsymb} at 27.75 0  28.25 0 /
  \plot 27.75 0  28.15 0 /
  
  \multiput{$\smallsq5$} at 31.5 -.5  31.5 0  31.5 .5  32 -1.5  32 -1  32 -.5 
  32 0  32 .5  32 1 /
  \multiput{\bulletsymb} at 31.75 -.5  32.25 -.5 /
  \plot 31.75 -.5  32.15 -.5 /
  
  \multiput{$\smallsq5$} at 35.75 -1  35.75 -.5  35.75 0  35.75 .5  35.75 1 /
  \multiput{\bulletsymb} at 36 -.5 /
  \arr{24.75 -1.5} {25.5 -3}
  \arr{28.75 -1.5} {29.5 -3}
  \arr{32.75 -1.5} {33.5 -3}
  \arr{26.5 -3} {27.25 -1.5}
  \arr{30.5 -3} {31.25 -1.5}
  \arr{34.5 -3} {35.25 -1.5}
  \put{$\scriptstyle(2)$} at 28 -4.75
  \put{$\scriptstyle(3)$} at 32 -4.75
  \setdots<2pt>
  \plot 26.5 -4.25  29.5 -4.25 /
  \plot 30.5 -4.25  33.5 -4.25 /
  \plot 20.5 -.25   23.5 -.25 /
  \plot 36.5 -.25   39.5 -.25 /
  \multiput{$\cdots$} at 22 3  26 3  30 3  34 3  38 3 / 
  \endpicture$$
  The number in parantheses describes if in the Auslander-Reiten sequence labelled,
  (1) the sequence of submodules, (2) the sequence of big modules, or (3) the
  sequence of the factors is not split exact.

\end{example}

\section{Isomorphy of Auslander-Reiten quivers.}  \label{section-isomorphy}

We demonstrate the result in the introduction about the independence
of the Auslander-Reiten quiver from the commutative uniserial base ring. 
As an immediate application we obtain the Auslander-Reiten quiver for the category
in Example~\ref{example-chains} (chains of subgroups). 
We use the terminology related to translation quivers and coverings from 
\cite{bongartz-gabriel}. 

\begin{definition}
  A full connected finite nonempty
  subquiver $\mathcal L$ of a (valued) translation quiver $\Gamma$ is a {\it slice}
  provided whenever $x\in\mathcal L$ is a noninjective point 
  and $f:x\to y$ an arrow 
  into a nonprojective point $y$, then exactly one of the arrows 
  $f$ or $\sigma(f):\tau y\to x$ is in $\mathcal L$. Here, $\sigma$ denotes the
  semitranslation in $\Gamma$.
\end{definition}

\begin{proposition} \label{covering}
  Let $\Lambda$ be a commutative local uniserial ring, $\mathcal P$ a finite
  poset, $\Gamma$ a connected component of the Auslander-Reiten quiver
  of $\sub_\Lambda\mathcal P$, and $\mathcal L$ a slice in $\Gamma$.
  Let $\tilde\Gamma\sto\pi\Gamma$ be the universal covering for $\Gamma$ and
  define the {\it type} of a point $x$ in $\tilde\Gamma$ as the type of
  $\pi(x)$ in $\Gamma$. 
  The structure of the (valued) translation quiver $\tilde\Gamma$ and the types 
  assigned to its points are uniquely determined by the structure of  $\mathcal L$
  as a (valued) quiver and by the types of the points in $\mathcal L$.
\end{proposition}

\begin{proof}
  A point $x\in\tilde\Gamma$ such that $\pi(x)\in\mathcal L$ defines a unique subquiver
  $\tilde{\mathcal L}\subset \tilde\Gamma$ which corresponds to $\mathcal L$ under $\pi$;
  this quiver is a slice in $\tilde\Gamma$. 
  We define a sequence $\mathcal L_i$ of slices shifted agains $\tilde{\mathcal L}$ 
  such that every point $y\in\tilde\Gamma$ with a path $y\to \tilde{\mathcal L}$ 
  (or a path $\tilde{\mathcal L}\to y$) is contained in one of the $\mathcal L_i$, $i\geq 0$
  (or in one of the $\mathcal L_i$, $i\leq 0$).
  We show that the types of the points in each of the $\mathcal L_i$, $i>0$, are determined uniquely by
  the types of the points in $\tilde{\mathcal L}$ and omit the proof of the corresponding result 
  for $\mathcal L_i$, $i<0$. 
  Given $i\geq 0$ and $z\in L_i$ a sink, obtain $L_{i+1}$ by performing the
  first possible of the following operations.
  a) First we deal with the case that $z$ is a projective vertex in $\tilde \Gamma$.  
    If $z$ is an
    isolated point then $\Lambda$ is a field, $z=[(\Lambda;(0)_{j\in\mathcal P})]$, and 
    we are done with all points $x$ which have a path into $\tilde{\mathcal L}$. 
    Otherwise there is a 
    unique arrow $y\to z$ in $\tilde \Gamma$ ending in $z$ 
    where $y$ is the class of the radical of the projective object
    corresponding to $z$ as noted in Observation~\ref{observation}. 
    This arrow is in $\mathcal L_i$ since $\mathcal L_i$ is connected.
    Let $\mathcal L_{i+1}$ be obtained from $\mathcal L_i$ by deleting the point $z$ and the 
    arrow $y\to z$. The graph $\mathcal L_{i+1}$ is connected, hence a slice. 

  b) Here we assume that $z$ is the 
    end term of a source map $y\to z$ where $y$ is some injective object 
    not in $\mathcal L_i$.  
    We have seen in Observation~\ref{observation} that the endterm $z$ 
    of the source map for $y$ is 
    an indecomposable object which is determined uniquely by its type.
    By adding to $\mathcal L_i$ the point $y$ and the arrow $y\to z$ 
    we arrive at a slice $\mathcal L_{i+1}$ in which the types of all points are determined 
    uniquely by the types of the points in $\mathcal L_i$. 

  c) If neither of the above two operations can be performed then we are in the situation that
    $z$ is the class of an endterm of an Auslander-Reiten sequence $\mathcal E$ 
    and that for each indecomposable
    summand $y$ of the middle term there is an arrow $y\to z$ in $\mathcal L_i$. 
    In this case we replace $z$ by $\tau z$ and each arrow $y\to z$ by a corresponding arrow 
    $\tau z\to y$ to obtain a new slice $\mathcal L_{i+1}$.  
    The type of $z$ decides whether the
    sequence $\mathcal E$ is split exact in each component (in which case additivity of types 
    yields the type of $\tau z$) or if $\mathcal E$ is 
    one of the sequences studied in the previous
    section (with the type of $\tau z$ known).
\end{proof}

The fact used in the above proof namely that the first term in the 
  source map for a projective module, and the last term in the sink map for an injective
  module are indecomposable modules (or zero), has the following consequence for the
shape of the Auslander-Reiten quiver of a submodule category.

\begin{proposition} \label{structure}
  Let $\Lambda$ be a commutative uniserial ring, $\mathcal P$ a finite
  poset and $\Gamma^s$ the stable part of a connected component $\Gamma$ 
  of the Auslander-Reiten quiver of $\sub_\Lambda\mathcal P$.
  If $\Gamma^s$ is non empty,
  then $\Gamma^s$ itself is connected and $\Gamma$ is obtained from $\Gamma^s$ by
  successively attaching non stable orbits. \qed
\end{proposition}

We can now give the proof of Theorem~\ref{introthm}.

\begin{proof}
We use the notation in the theorem.  Since the slices $\mathcal L_\Lambda$ and
$\mathcal L_\Delta$ are isomorphic as quivers and consist of objects of the same types,
Proposition~\ref{covering} yields  universal coverings 
$\pi_\Lambda:\tilde\Gamma_\Lambda\to \Gamma_\Lambda$,
$\pi_\Delta:\tilde\Gamma_\Delta\to \Gamma_\Delta$ and a type preserving isomorphism
of translation quivers $\tilde\varphi$ which maps $\mathcal L_\Lambda$ onto $\mathcal L_\Delta$.
$$\begin{CD} \tilde\Gamma_\Lambda @>\tilde\varphi>> \tilde\Gamma_\Delta \\
  @V\pi_\Lambda VV   @VV\pi_\Delta V \\ \Gamma_\Lambda @. \Gamma_\Delta \end{CD}$$
We show that there is a type preserving isomorphism $\varphi:\Gamma_\Lambda\to\Gamma_\Delta$
of translation quivers which makes the diagram commutative. 

\smallskip
Suppose that $x_\Lambda\in\Gamma_\Lambda$ and $x_\Delta\in\Gamma_\Delta$ are objects 
which have the same type and which are determined uniquely by their type. 
For example, projective objects, injective objects, and the objects in $\mathcal L$ 
have this property.
Then the fibers $\pi_\Lambda^{-1}(x_\Lambda)$ and $\pi_\Delta^{-1}(x_\Delta)$
coincide and we can define $\varphi(x_\Lambda)=x_\Delta$. 
Since $\pi_\Lambda$, $\pi_\Delta$ define coverings, the maps commute with the translation
and the definition of $\varphi$ can be extended to the preprojective, to the preinjective, and
to the stable objects. Thus, $\varphi$ is a bijection
of the points in $\Gamma_\Lambda$ and $\Gamma_\Delta$ which satisfies 
$\varphi\circ\pi_\Lambda=\pi_\Delta\circ\tilde\varphi$ and which 
extends to an isomorphism of translation quivers.
\end{proof}

\begin{remark}
  Note that we cannot expect in Theorem~\ref{introthm} that
  the type determines an object uniquely,
  up to isomorphism, not even in the representation finite case,
  as the following example shows.
  Thus, condition (3) may not hold for all objects in the component
  $\mathcal C$.
\end{remark}

\begin{example} \label{type-example}
  Let $\mathcal P=\bullet_1$   and $\Lambda=\Zmod5$.
The category $\sub_\Lambda\mathcal P$ has finite representation type, in fact all 50
indecomposables have been determined in \cite{richman}.
The Auslander-Reiten quiver $\Gamma$ has a slice which consists of points which
are determined uniquely by their type.  But not all points in  $\Gamma$ are determined 
by their type:   The two indecomposable   objects in $\sub_\Lambda\mathcal P$,
  \begin{eqnarray*}
    (M_*;M_1)& =&\left(\Zmod5\oplus\Zmod2;\;\Zmod3([p^2],[1])\right)
    \quad \mbox{and}\\
    (M'_*;M'_1) &= &\left(\Zmod5\oplus\Zmod2;\;\Zmod3([p^2],[p])\right),
  \end{eqnarray*}
  are pictured as follows.
  $$\beginpicture 
  \setcoordinatesystem units <0.3cm,0.3cm>
  \put{$(M_*;M_1):$} at -4 2.5
  \multiput{\sq} at 0 5  0 4  0 3  0 2  0 1  1 3  1 2 /
  \put{$\scriptstyle \bullet$} at 0.5 3 
  \put{$\scriptstyle \bullet$} at 1.5 3 
  \plot 0.5 3  1.5 3 /
  \endpicture\qquad\qquad
  \beginpicture 
  \setcoordinatesystem units <0.3cm,0.3cm>
  \put{$(M'_*;M'_1):$} at -4 2.5
  \multiput{\sq} at 0 5  0 4  0 3  0 2  0 1  1 4  1 3 /
  \put{$\scriptstyle \bullet$} at 0.5 3 
  \put{$\scriptstyle \bullet$} at 1.5 3 
  \plot 0.5 3  1.5 3 /
  \endpicture$$
  Both objects have type $(52;3)$ but since their quotients
  $M_*/M_1\cong\Zmod4$ and
  $M_*'/M_1'\cong\Zmod3\oplus\Zmod1$ are not isomorphic in $\mod\Lambda$, 
  the objects $M$ and $M'$ cannot   be isomorphic in $\sub_\Lambda\mathcal P$.
\end{example}

\begin{example}[Chains of subgroups, III] We continue Example~\ref{chains-two} where we
  have seen that the   sequence of irreducible maps 
  $\vcenter{\beginpicture 
    \setcoordinatesystem units <.7cm,.7cm>
    \put{$\smallsq5$} at 0 0
    \multiput{} at 0 -.5  .5 .5 /
    \endpicture}\to 
  \vcenter{\beginpicture 
    \setcoordinatesystem units <.7cm,.7cm>
    \put{$\smallsq5$} at 0 0
    \multiput{} at 0 -.5  .5 .5 /
    \put{$\eyesm$} at .25 0
    \endpicture}\to
  \vcenter{\beginpicture 
    \setcoordinatesystem units <.7cm,.7cm>
    \put{$\smallsq5$} at 0 0
    \multiput{} at 0 -.5  .5 .5 /
    \put{$\eyesm$} at .25 0
    \put{$\nosem$} at .25 0
    \endpicture}\to
  \vcenter{\beginpicture 
    \setcoordinatesystem units <.7cm,.7cm>
    \put{$\smallsq5$} at 0 0
    \multiput{} at 0 -.5  .5 .5 /
    \put{$\eyesm$} at .25 0
    \put{$\nosem$} at .25 0
    \put{$\mouthm$} at .25 0
    \endpicture}$
  forms a slice in the Auslander-Reiten quiver for $\sub_\Lambda\mathcal P$.
  Moreover, as endterms of Auslander-Reiten sequences   which are not split exact
  in each component, the modules are uniquely determined by their type.  
  In fact, every indecomposable object is determined uniquely by its type and the
  computation of the Auslander-Reiten quiver from the slice is a straightforward task.
  The example serves as an illustration for Theorem~\ref{introthm} as for any two
  commutative uniserial rings $\Lambda$ and $\Delta$ of length two, there is a 
  type preserving isomorphism of translation quivers between the 
  Auslander-Reiten quivers for $\sub_\Lambda\mathcal P$ and 
  $\sub_\Delta\mathcal P$.  Each quiver is as follows.
  $$
  \vcenter{\beginpicture
    \setcoordinatesystem units <.7cm,.7cm> 
    \put{} at  0 0
    \put{} at  15 5
    \multiput{$\oneo$} at 1 1  2 2  3 3  4 4  6 4  7 3  8 2  9 1  11 1  12 2 
    13 3  14 4  /
    \multiput{$\oo$} at 0 4  1 5  0 2  1 3  2 4  3 5  2 0  3 1  4 2  5 3  7 5  
    4 0  5 1  6 2  8 4  9 5  6 0  7 1  9 3  10 4  8 0  10 2  
    11 5  10 2  11 3  12 4  13 5  12 0  13 1  14 2  15 3  14 0  15 1 /
    \multiput{\eyesm} at 2 2  3 3  4 4  7 3  8 2  9 1  12 2  13 3  14 4 /
    \multiput{\eyesl} at 3 1  4 2  5 3  8 4  9 3  10 2  0 2  13 1  14 2  15 3 /
    \multiput{\eyesh} at 4 0  5 1  6 2  6 0  7 1  8 0  9 5  10 4  0 4  1 5  1 3  
    2 4  3 5  14 0  15 1  11 5  11 3  12 4  13 5  /
    \multiput{\nosem} at 3 3  4 4  8 2  9 1  13 3  14 4  /
    \multiput{\nosel} at 4 2  5 3  9 3  10 2  0 2  0 4  1 3  5 1  6 2  10 4 
    14 2  15 3  10 2  11 3  15 1 /
    \multiput{\noseh} at 7 1  8 0  1 5  2 4  3 5  6 0  11 5  12 4  13 5 /
    \multiput{\mouthm} at 4 4  9 1  14 4 /
    \multiput{\mouthl} at 5 3  10 2  0 2  1 3  2 4  6 2  7 1  15 3  11 3  12 4 /
    \multiput{\mouthh} at 8 0  3 5  13 5  /
    \arr{0.35 4.35} {0.65 4.65}
    \arr{0.35 2.35} {0.65 2.65}
    \arr{1.35 3.35} {1.65 3.65}
    \arr{2.35 4.35} {2.65 4.65}
    \arr{1.35 1.35} {1.65 1.65}
    \arr{2.35 2.35} {2.65 2.65}
    \arr{3.35 3.35} {3.65 3.65}
    \arr{2.35 0.35} {2.65 0.65}
    \arr{3.35 1.35} {3.65 1.65}
    \arr{4.35 2.35} {4.65 2.65}
    \arr{4.35 0.35} {4.65 0.65}
    \arr{5.35 1.35} {5.65 1.65}
    \arr{6.35 2.35} {6.65 2.65}
    \arr{5.35 3.35} {5.65 3.65}
    \arr{6.35 4.35} {6.65 4.65}
    \arr{7.35 3.35} {7.65 3.65}
    \arr{8.35 4.35} {8.65 4.65}
    \arr{6.35 0.35} {6.65 0.65}
    \arr{7.35 1.35} {7.65 1.65}
    \arr{8.35 2.35} {8.65 2.65}
    \arr{9.35 3.35} {9.65 3.65}
    \arr{8.35 0.35} {8.65 0.65}
    \arr{9.35 1.35} {9.65 1.65}
    \arr{10.35 4.35} {10.65 4.65}
    \arr{10.35 2.35} {10.65 2.65}
    \arr{11.35 3.35} {11.65 3.65}
    \arr{12.35 4.35} {12.65 4.65}
    \arr{11.35 1.35} {11.65 1.65}
    \arr{12.35 2.35} {12.65 2.65}
    \arr{13.35 3.35} {13.65 3.65}
    \arr{12.35 0.35} {12.65 0.65}
    \arr{13.35 1.35} {13.65 1.65}
    \arr{14.35 2.35} {14.65 2.65}
    \arr{14.35 0.35} {14.65 0.65}
    
    \arr{0.35 3.65} {0.65 3.35}
    \arr{1.35 4.65} {1.65 4.35}
    \arr{0.35 1.65} {0.65 1.35}
    \arr{1.35 2.65} {1.65 2.35}
    \arr{2.35 3.65} {2.65 3.35}
    \arr{3.35 4.65} {3.65 4.35}
    \arr{1.35 0.65} {1.65 0.35}
    \arr{2.35 1.65} {2.65 1.35}
    \arr{3.35 2.65} {3.65 2.35}
    \arr{4.35 3.65} {4.65 3.35}
    \arr{3.35 0.65} {3.65 0.35}
    \arr{4.35 1.65} {4.65 1.35}
    \arr{5.35 2.65} {5.65 2.35}
    \arr{6.35 3.65} {6.65 3.35}
    \arr{7.35 4.65} {7.65 4.35}
    \arr{5.35 0.65} {5.65 0.35}
    \arr{6.35 1.65} {6.65 1.35}
    \arr{7.35 2.65} {7.65 2.35}
    \arr{8.35 3.65} {8.65 3.35}
    \arr{9.35 4.65} {9.65 4.35}
    \arr{7.35 0.65} {7.65 0.35}
    \arr{8.35 1.65} {8.65 1.35}
    \arr{9.35 2.65} {9.65 2.35}
    \arr{10.35 3.65} {10.65 3.35}
    \arr{11.35 4.65} {11.65 4.35}
    \arr{10.35 1.65} {10.65 1.35}
    \arr{11.35 2.65} {11.65 2.35}
    \arr{12.35 3.65} {12.65 3.35}
    \arr{13.35 4.65} {13.65 4.35}
    \arr{11.35 0.65} {11.65 0.35}
    \arr{12.35 1.65} {12.65 1.35}
    \arr{13.35 2.65} {13.65 2.35}
    \arr{14.35 3.65} {14.65 3.35}
    \arr{13.35 0.65} {13.65 0.35}
    \arr{14.35 1.65} {14.65 1.35}

    \setdots<2pt>
    \plot -.5 1  .5 1 /
    \plot 4.5 4  5.5 4 /
    \plot 9.5 1  10.5 1 /
    \plot 14.5 4  15.5 4 /
    
    \multiput{$\cdots$} at 15.8 1  15.8 3   -.8 2  -.8 4 /
    \endpicture}
  $$
\end{example}

\section{Submodule categories of type $\mathcal S_3(\Lambda)$
where $\Lambda$ has length 6.}\label{section-test}

For $\Lambda$ a commutative uniserial ring of length 6 we show that the category 
$\mathcal S_3(\Lambda)$ is representation finite. 
We compute its Auslander-Reiten quiver and demonstrate that this quiver is independent
of the choice of $\Lambda$.  In particular, there are no parametrized families
of indecomposable objects in $\mathcal S_3(\mathbb Z/p^6)$.
This finishes our demonstration that the first occurances of parametrized
families of subgroup embeddings are in the categories $\mathcal S_4(\mathbb Z/p^6)$
and $\mathcal S_3(\mathbb Z/p^7)$. 

\smallskip
We proceed as follows.
The corresponding category $\mathcal S_3(k[x]/x^6)$ where $k=\Lambda/\rad \Lambda$ 
has 84 indecomposables
and the Auslander-Reiten quiver is as pictured on page~\pageref{arq36}.
We have seen in Example~\ref{bounded-I} that the modules in the top orbit and in the non stable 
orbit occur as endterms of Auslander-Reiten sequences which are not split exact
in each component; thus, also the corresponding modules in the category 
$\mathcal S_3(\Lambda)$ form orbits under the Auslander-Reiten translation.  
Starting from these modules, we construct a slice in the Auslander-Reiten quiver for
$\mathcal S_3(\Lambda)$; as a quiver, the slice is isomorphic to a corresponding 
slice for $\mathcal S_3(k[x]/x^6)$.  We deduce that there is a type preserving 
isomorphism between the components of the two Auslander-Reiten
quivers.  Thus we have detected a finite component, and the Harada-Sai lemma implies that 
we have computed  the complete Auslander-Reiten quiver for 
$\mathcal S_3(\Lambda)$.

\medskip
First we adapt Theorem~\ref{introthm} to the case of bounded subgroups.

\begin{corollary}
  Suppose $\Lambda,\Delta$ are commutative
  uniserial rings of the same length $n$, $m\leq n$ is a natural number, 
  $\Gamma_\Lambda$ and $\Gamma_\Delta$ are connected components 
  of the Auslander-Reiten quivers of $\mathcal S_m(\Lambda)$ and
  $\mathcal S_m(\Delta)$, and $\mathcal L_\Lambda$ and $\mathcal L_\Delta$
  are slices in $\Gamma_\Lambda$ and $\Gamma_\Delta$, respectively, such that the following 
  conditions are satisfied:
  \begin{enumerate}
  \item The slices $\mathcal L_\Lambda$ and $\mathcal L_\Delta$ are isomorphic as 
    valued graphs.
  \item Points in $\mathcal L_\Lambda$ and $\mathcal L_\Delta$ which correspond to each other
    under this isomorphism represent indecomposable objects of the same type.
  \item Each indecomposable object represented by a point in $\mathcal L_\Lambda$ or in
    $\mathcal L_\Delta$ is determined uniquely, up to isomorphism, by its type.
  \end{enumerate}
  Then the components $\Gamma_\Lambda$ and $\Gamma_\Delta$ are isomorphic
  as valued translation quivers, 
  and points which correspond to each other under this isomorphism
  represent indecomposable objects of the same type.
\end{corollary}

\begin{proof}
  The proof of Theorem~\ref{introthm} in Section~\ref{section-isomorphy} applies also 
  in this situation:  In \cite{ms-bounded} the projective and the injective indecomposables 
  in categories of type $\mathcal S_m(\Lambda)$ and their respective sink and source maps
  have been computed.  As in the case of submodule representations of a poset, each source 
  or target in such a map is zero or an indecomposable module which is determined uniquely
  by its type.  The second fact needed in the proof is the structure of those Auslander-Reiten
  sequences which are not split exact in each component.  These sequences have been 
  constructed in Example~\ref{bounded-I}, (1) and (2). 
\end{proof}

\medskip
Next we provide us with a tool to verify that given short exact sequences are 
in fact Auslander-Reiten sequences. 
Recall from \cite[Proposition V.2.2]{ars} that in a category of modules over
an artin algebra, a nonsplit short exact sequence
$0\to D\text{Tr}\,C\sto fB\sto gC\to 0$ where $C$ is an indecomposable
nonprojective module is an Auslander-Reiten
sequence provided only that every endomorphism of $C$ which is not
an automorphism factors factors over $g$. We adapt this 
result to our situation where we are dealing with a full subcategory
$\mathcal S$ of a module category which is closed under extensions.  

\smallskip
We first recall the result corresponding \cite[Proposition V.2.1]{ars}.

\begin{lemma}
  Let $C$ be an indecomposable object in ${\mathcal S}$ which
  is not $\Ext$-projective, and let 
  $0\to A\to B\to C\to 0$ be an Auslander-Reiten sequence 
  in ${\mathcal S}$.
  Then $\Ext^1_R(C,A)$ has a simple socle, both as 
  $\End(C)$- and as $\End(A)$-module.  
  These socles coincide and each nonzero element
  in the socle is an Auslander-Reiten sequence 
  in ${\mathcal S}$. \qed \end{lemma}

\medskip
The following result is a minor adaption of 
\cite[Proposition V.2.2]{ars} to the case of full subcategories.

\begin{proposition} \label{ar-test}
  Let $\mathcal S$ be a full extension closed subcategory of a category of
  modules over an artin algebra.
  For $C$ an indecomposable object in ${\mathcal S}$ which
  is not $\Ext$-projective, let 
  $0\to A\to B'\to C\to 0$ be an Auslander-Reiten sequence 
  in ${\mathcal S}$.  The following assertions are equivalent for
  a nonsplit short exact sequence 
  $${\mathcal E}: \;0\lto{}A\lto{f}B\lto{g}C\lto{}0.$$
  \begin{enumerate}
  \item[$1.$] The sequence ${\mathcal E}$ is an 
    Auslander-Reiten sequence in ${\mathcal S}$.
  \item[$2.$] Every nonautomorphism $h:C\to C$ factors
    through $C$. 
  \item[$3.$] The image of the map $\Hom(C,g):\Hom(C,B)\to 
    \Hom(C,C)$ is $\Rad\End(C)$. 
  \item[$4.$] The contravariant defect 
    ${\mathcal E}^*(C) = \End(C)/\Im\Hom(C,g)$
    is a simple $\End(C)$-module.
  \item[$5.$] The socle of the $\End(C)$-module $\Ext^1_R(C,A)$ is generated by the
    class of $\mathcal E$.    
  \item[$2'.$] Every nonautomorphism $k:A\to A$ factors through $f$.
  \item[$3'.$] The the map $\Hom(f,A):\Hom(B,A)\to \Hom(A,A)$
    has image $\Rad\End(A)$. 
  \item[$4'.$] The covariant defect 
    ${\mathcal E}_*(A) = \End(C)/\Im\Hom(f,A)$ is a simple 
    $\End(A)$-module.
  \item[$5'.$] The socle of the $\End(A)$-module $\Ext^1_R(C,A)$ is generated by 
    the class of $\mathcal E$.
  \end{enumerate}
\end{proposition}

\begin{proof}
  The implications $1.\Rightarrow 2.$ and $2.\Leftrightarrow3.\Leftrightarrow4.$
  are obvious and $1.\Leftrightarrow5.$ follows from Lemma 4.1.
  We only show that $4.\Rightarrow5.$ 
  Consider the long exact sequence
  $$0\to \Hom(C,A)\to \Hom(C,B)\to \Hom(C,C)\sto\delta\Ext^1_R(C,A)\to\cdots$$
  Recall that $\delta$ is given by sending the identity map $1_C$ to
  the class of $\mathcal E$ in $\Ext^1_R(C,A)$ and that the image of $\delta$
  is the contravariant
  defect ${\mathcal E}^*(C)$.
  Thus, ${\mathcal E}^*(C)$ is the submodule of $\Ext_R^1(C,A)$
  generated by $\bar{\mathcal E}$. If ${\mathcal E}^*(C)_{\End(C)}$ is a simple module
  then $\bar{\mathcal E}$ is a socle generator by Lemma 4.1.
  The statements $2'.$, $3'.$, $4'.$ and $5'.$ follow by duality.
\end{proof}

\medskip 
We return to the setup from Example~\ref{bounded-I} 
and compute the Auslander-Reiten quiver for 
$\mathcal S_3(\mathbb Z/p^6)$.

\begin{theorem}
  Let $\Lambda$ be a commutative uniserial ring of length 6.
  \begin{enumerate}
  \item The category $\mathcal S_3(\Lambda)$ has 84 
    isomorphism classes of indecomposable objects. The stable part of the Auslander-Reiten quiver
    has type $\mathbb Z\mathbb E_8/\tau^{10}$ and there are two nonstable orbits of length
    1 and 3 attached to the boundaries of the wings of width 5 and 3, respectively.
  \item Each indecomposable object in $\mathcal S_3(\Lambda)$ is determined
    uniquely, up to isomorphism, by its type with the following 10 exceptions
    for which we specify the type and the cotype.
    $$
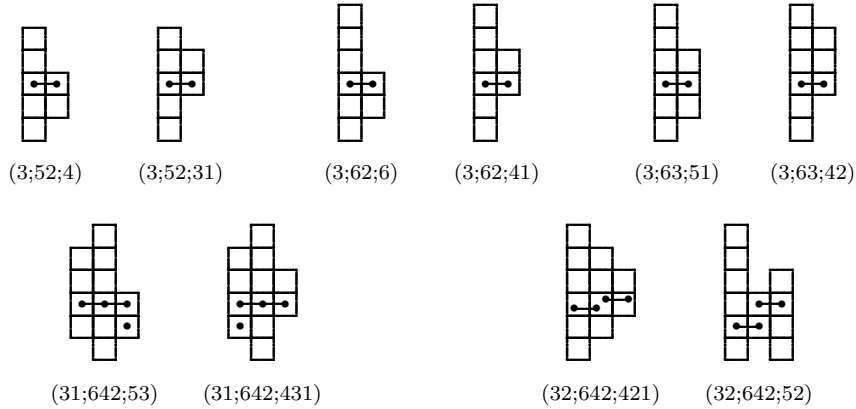
\beginpicture\setcoordinatesystem units <3mm,3mm>
    \multiput{$\sq$} at 2 10  2 11  2 12  2 13  2 14  3 11  3 12 /
    \multiput{$\scriptscriptstyle\bullet$} at 2.5 12  3.5 12 /
    \plot 2.5 12  3.5 12 /
    \put{$\scriptstyle(3;52;4)$} at 3 8
    \multiput{$\sq$} at 8 10  8 11  8 12  8 13  8 14  9 12  9 13 /
    \multiput{$\scriptscriptstyle\bullet$} at 8.5 12  9.5 12 /
    \plot 8.5 12  9.5 12 /
    \put{$\scriptstyle(3;52;31)$} at 9 8
    \multiput{$\sq$} at 16 10  16 11  16 12  16 13  16 14  16 15  17 11  17 12 /
    \multiput{$\scriptscriptstyle\bullet$} at 16.5 12  17.5 12 /
    \plot 16.5 12  17.5 12 /
    \put{$\scriptstyle(3;62;6)$} at 17 8
    \multiput{$\sq$} at 22 10  22 11  22 12  22 13  22 14  22 15  23 12  23 13 /
    \multiput{$\scriptscriptstyle\bullet$} at 22.5 12  23.5 12 /
    \plot 22.5 12  23.5 12 /
    \put{$\scriptstyle(3;62;41)$} at 23 8
    \multiput{$\sq$} at 30 10  30 11  30 12  30 13  30 14  30 15  31 11  31 12  31 13 /
    \multiput{$\scriptscriptstyle\bullet$} at 30.5 12  31.5 12 /
    \plot 30.5 12  31.5 12 /
    \put{$\scriptstyle(3;63;51)$} at 31 8
    \multiput{$\sq$} at 36 10  36 11  36 12  36 13  36 14  36 15  37 12  37 13  37 14 /
    \multiput{$\scriptscriptstyle\bullet$} at 36.5 12  37.5 12 /
    \plot 36.5 12  37.5 12 /
    \put{$\scriptstyle(3;63;42)$} at 37 8
    \endpicture$$
    $$\beginpicture\setcoordinatesystem units <3mm,3mm>
    \multiput{$\sq$} at 4 1  4 2  4 3  4 4  5 0  5 1  5 2  5 3  5 4  5 5  6 1  6 2 /
    \multiput{$\scriptscriptstyle\bullet$} at 4.5 2  5.5 2  6.5 2  6.5 1 /
    \plot 4.5 2  6.5 2 /
    \put{$\scriptstyle(31;642;53)$} at 5.5 -2
    \multiput{$\sq$} at 11 1  11 2  11 3  11 4  12 0  12 1  12 2  12 3  12 4  
    12 5  13 2  13 3 /
    \multiput{$\scriptscriptstyle\bullet$} at 11.5 2  12.5 2  13.5 2  11.5 1 /
    \plot 11.5 2  13.5 2 /
    \put{$\scriptstyle(31;642;431)$} at 12.5 -2
    \multiput{$\sq$} at 26 0  26 1  26 2  26 3  26 4  26 5  
    27 1  27 2  27 3  27 4  28 2  28 3  /
    \multiput{$\scriptscriptstyle\bullet$} at 26.3 1.8  27.3 1.8  27.7 2.2  28.7 2.2 /
    \plot 26.3 1.8  27.3 1.8 /
    \plot 27.7 2.2  28.7 2.2 /
    \put{$\scriptstyle(32;642;421)$} at 27.5 -2
    \multiput{$\sq$} at 33 0  33 1  33 2  33 3  33 4  33 5 
    34 1  34 2  35 0  35 1  35 2  35 3  /
    \multiput{$\scriptscriptstyle\bullet$} at 33.5 1  34.5 1  34.5 2  35.5 2 /
    \plot 33.5 1  34.5 1 /
    \plot 34.5 2  35.5 2 /
    \put{$\scriptstyle(32;642;52)$} at 34.5 -2
    \endpicture$$
  \end{enumerate}
\end{theorem}

\medskip
\begin{proof}
  Let $k$ be any field.  We show that there is a type preserving isomorphism of 
  (valued) translation quivers between the Auslander-Reiten quiver for
  $\mathcal S_3(\Lambda)$ and the Auslander-Reiten quiver for $\mathcal S_3(k[x]/x^6)$,
  of which we include a copy on page~\pageref{arq36}.
  Both assertions in the theorem follow easily.

  \smallskip
  The diagram below is part of the Auslander-Reiten quiver for 
  the category $\mathcal S_3(k[x]/x^6)$
  and has been determined in \cite{ms-bounded} using covering theory. 
  We show that the shaded region defines a slice also in the
  Auslander-Reiten quiver for $\mathcal S_3(\Lambda)$.
\begin{figure}[ht] 
  $$\beginpicture\setcoordinatesystem units <4mm,4mm>
  \multiput{} at 14 24  44 -5 /
  \put{A slice in the\strut} at 40 6.5
  \put{Auslander-Reiten\strut} at 40 5.25
  \put{quiver for $\mathcal S_3(\Lambda)$\strut} at 40 4
\multiput{$\smallsq3$} at  25.85 -4.6  25.85 -4.3  25.85 -4 /
\multiput{\bulletsymb} at 26 -4 /

\multiput{$\smallsq3$} at 29.85 -4.3  29.85 -4  29.85 -3.7 /

\multiput{$\smallsq3$} at 33.85 -4.9  33.85 -4.6  33.85 -4.3  33.85 -4
33.85 -3.7  33.85 -3.4 /
\multiput{\bulletsymb} at 34 -4.3 /
\multiput{$\smallsq3$} at 23.85 -.6  23.85 -.3  23.85 0  /
\multiput{\bulletsymb} at 24 -.3 /

\multiput{$\smallsq3$} at 27.7 -.6  27.7 -.3  27.7 0  27.7 .3 28 -.3  28 0 /
\multiput{\bulletsymb} at 27.85 0  28.15 0 /
\plot 27.85 0  28.15 0 /

\multiput{$\smallsq3$} at 31.7 -.3  31.7 0  31.7 .3  32 -.9  32 -.6  32 -.3 
32 0  32 .3  32 .6 /
\multiput{\bulletsymb} at 31.85 -.3  32.15 -.3 /
\plot 31.85 -.3  32.15 -.3 /
\multiput{$\smallsq3$} at 25.7 3.4  25.7 3.7  25.7 4  25.7 4.3  26 3.7  26 4 /
\multiput{\bulletsymb} at 25.85 3.7  26.15 3.7 /
\plot  25.85 3.7  26.15 3.7 /

\multiput{$\smallsq3$} at 29.55 3.4  29.55 3.7  29.55 4  29.55 4.3  29.85 3.7
29.85 4  30.15 3.1  30.15 3.4  30.15 3.7  30.15 4  30.15 4.3
30.15 4.6 /
\multiput{\bulletsymb} at 29.7 4  30 4  30 3.7  30.3 3.7 /
\plot 29.7 4  30 4 /
\plot  30 3.7  30.3 3.7 /
  \multiput{$\smallsq3$} at 23.4 7.7  23.4 8  23.7 7.1  23.7 7.4  23.7 7.7
  23.7 8  23.7 8.3  23.7 8.6  24 7.4  24 7.7  24 8  24 8.3 
  24.3 7.7 /
  \multiput{\bulletsymb} at 23.5 7.65  23.8 7.65  24.1 7.65  24.2 7.75  
  24.5 7.75 /
  \plot 23.5 7.65  24.1 7.65 /
  \plot 24.2 7.75  24.5 7.75 /
  
\multiput{$\smallsq3$} at 27.55 7.1  27.55 7.4  27.55 7.7  27.55 8  27.55 8.3
27.55 8.6  27.85 7.4  27.85 7.7  27.85 8  27.85 8.3  28.15 7.7  28.15 8         /
\multiput{\bulletsymb} at 27.65 7.65  27.95 7.65  28.05 7.75  28.35 7.75 /
\plot 27.65 7.65  27.95 7.65 /
\plot   28.05 7.75  28.35 7.75 /

\multiput{$\smallsq3$} at 31.55 7.4  31.55 7.7  31.55 8  31.55 8.3  31.55 8.6
31.85 7.7  31.85 8  32.15 7.4  32.15 7.7  32.15 8  32.15 8.3 
/
\multiput{\bulletsymb} at 31.7 7.7  32 7.7  32 8 32.3 8 /
\plot 31.7 7.7  32 7.7 /
\plot  32 8 32.3 8 /
  \multiput{$\smallsq3$} at 25.7 6.1  25.7 6.4  25.7 6.7  25.7 7  25.7 7.3 
  25.7 7.6  26 6.7  26 7 /
  \multiput{\bulletsymb} at 25.85 6.7  26.15 6.7 /
  \plot  25.85 6.7  26.15 6.7 /
  
\multiput{$\smallsq3$} at 29.85 6.4  29.85 6.7  29.85 7  29.85 7.3 /
\multiput{\bulletsymb} at 30 6.7 /
  \multiput{$\smallsq3$} at 21.55 10.7  21.55 11  21.85 10.1  21.85 10.4
  21.85 10.7  21.85 11  21.85 11.3  21.85 11.6  
  22.15 10.4  22.15 10.7  22.15 11  22.15 11.3 /
  \multiput{\bulletsymb} at 21.7 10.7  22 10.7  22.3 10.7  22.3 10.4 /
  \plot 21.7 10.7  22.3 10.7 /
  
\multiput{$\smallsq3$} at 25.55 10.1  25.55 10.4  25.55 10.7  25.55 11 
25.55 11.3  25.55 11.6  25.85 10.4  25.85 10.7  25.85 11
25.85 11.3  26.15 10.7 / 
\multiput{\bulletsymb} at 25.65 10.65  25.95 10.65  26.05 10.75  26.35 10.75 /
\plot 26.05 10.75  26.35 10.75 /
\plot 25.65 10.65  25.95 10.65 /
  
\multiput{$\smallsq3$} at 29.7 10.4  29.7 10.7  29.7 11 29.7 11.3  29.7 11.6 
30 10.7  30 11 /
\multiput{\bulletsymb} at 29.85 10.7  30.15 10.7 /
\plot 29.85 10.7  30.15 10.7 /

\multiput{$\smallsq3$} at 33.55 10.4  33.55 10.7  33.55 11  33.55 11.3  
33.55 11.6  33.55 11.9  33.85 10.7  33.85 11  34.15 10.4  34.15 10.7
34.15 11  34.15 11.3 /
\multiput{\bulletsymb} at 33.7 10.7  34 10.7  34 11  34.3 11 /
\plot 34 11  34.3 11 /
\plot  33.7 10.7  34 10.7 /
  \multiput{$\smallsq3$} at 19.55 13.7  19.55 14  19.85 13.1  19.85 13.4  
  19.85 13.7  19.85 14  19.85 14.3  19.85 14.6  20.15 13.4  20.15 13.7 
  20.15 14  20.15 14.3 /
  \multiput{\bulletsymb} at 19.7 13.7  20 13.7  20.3 13.7 /
  \plot 19.7 13.7  20.3 13.7 /
  
\multiput{$\smallsq3$} at 23.7 13.1  23.7 13.4  23.7 13.7  23.7 14 
23.7 14.3  23.7 14.6  24 13.4  24 13.7  24 14  24 14.3 /
\multiput{\bulletsymb} at 23.85 13.7  24.15 13.7  24.15 13.4 /
\plot 23.85 13.7  24.15 13.7 /

\multiput{$\smallsq3$} at 27.7 13.4  27.7 13.7  27.7 14  27.7 14.3  27.7 14.6
28 13.7 /
\multiput{\bulletsymb} at 27.85 13.7  28.15 13.7 /
\plot 27.85 13.7  28.15 13.7 /

\multiput{$\smallsq3$} at 31.7 13.4  31.7 13.7  31.7 14  31.7 14.3  31.7 14.6
31.7 14.9  32 13.7  32 14 /
\multiput{\bulletsymb} at 31.85 13.7  32.15 13.7 /
\plot 31.85 13.7  32.15 13.7 /

\multiput{$\smallsq3$} at 35.7 13.7  35.7 14  36 13.4  36 13.7  36 14  36 14.3 
/
\multiput{\bulletsymb} at 35.85 13.7  35.85 14  36.15 14 /
\plot 35.85 14  36.15 14 /
  \multiput{$\smallsq3$} at 17.7 16.7  17.7 17  18 16.1  18 16.4  18 16.7
  18 17  18 17.3 /
  \multiput{\bulletsymb} at 17.85 16.7  18.15 16.7 /
  \plot 17.85 16.7  18.15 16.7 /
  
\multiput{$\smallsq3$} at 21.7 16.1  21.7 16.4  21.7 16.7  21.7 17  21.7 17.3
21.7 17.6  22 16.4  22 16.7  22 17  22 17.3 /
\multiput{\bulletsymb} at 21.85 16.7  22.15 16.7 /
\plot 21.85 16.7  22.15 16.7 /

\multiput{$\smallsq3$} at 25.85 16.4  25.85 16.7  25.85 17  25.85 17.3
25.85 17.6 /
\multiput{\bulletsymb} at 26 16.4 /

\multiput{$\smallsq3$} at 29.7 16.4  29.7 16.7  29.7 17  29.7 17.3  29.7 17.6
29.7 17.9  30 16.7 /
\multiput{\bulletsymb} at 29.85 16.7  30.15 16.7 /
\plot  29.85 16.7  30.15 16.7 /

\multiput{$\smallsq3$} at 33.85 16.7  33.85 17 /
\multiput{\bulletsymb} at 34 16.7 /

\multiput{$\smallsq3$} at 37.7 17  38 16.4  38 16.7  38 17  38 17.3 /
\multiput{\bulletsymb} at 37.85 17  38.15 17 /
\plot 37.85 17  38.15 17 /

  \multiput{$\smallsq3$} at 15.85 19.7  15.85 20 /
  
\multiput{$\smallsq3$} at 19.85 19.1  19.85 19.4  19.85 19.7  19.85 20
19.85 20.3  /
\multiput{\bulletsymb} at 20 19.7 /

\multiput{$\smallsq3$} at 23.85 19.4  23.85 19.7  23.85 20  23.85 20.3
23.85 20.6 / 

\multiput{$\smallsq3$} at 27.85 19.4  27.85 19.7  27.85 20  27.85 20.3
27.85 20.6  27.85 20.9 /
\multiput{\bulletsymb} at 28 19.4 /

\multiput{$\smallsq3$} at 31.85 19.7 /
\multiput{\bulletsymb} at 32 19.7 /

\multiput{$\smallsq3$} at 35.85 20 /

\multiput{$\smallsq3$} at 39.85 19.4  39.85 19.7  39.85 20  39.85 20.3 /
\multiput{\bulletsymb} at 40 20 /
\multiput{$\smallsq3$} at 25.85 22.4  25.85 22.7  25.85 23  25.85 23.3
25.85 23.6  25.85 23.9 /
\arr{24.5 -1} {25.5 -3}
  \arr{28.5 -1} {29.5 -3}
  \arr{32.5 -1} {33.5 -3}
  \arr{26.5 -3} {27.5 -1}
  \arr{30.5 -3} {31.5 -1}
  \arr{34.5 -3} {35.5 -1}
  \arr{24.5 1} {25.25 2.5}
  \arr{28.5 1} {29.25 2.5}
  \arr{32.5 1} {33.25 2.5}
  \arr{22.75 2.5} {23.5 1}
  \arr{26.75 2.5} {27.5 1}
  \arr{30.75 2.5} {31.5 1}
  \arr{24.75 6.5} {25.5 5}
  \arr{28.75 6.5} {29.5 5}
  \arr{26.5 5} {27.25 6.5}
  \arr{30.5 5} {31.25 6.5}
  \arr{25 7.5} {25.5 7.25}
  \arr{29 7.5} {29.5 7.25}
  \arr{33 7.5} {34 7}
  \arr{22 7} {23 7.5}
  \arr{26.5 7.25} {27 7.5}
  \arr{30.5 7.25} {31 7.5}
  \arr{20.67 9} {21.33 10}
  \arr{24.67 9} {25.33 10}
  \arr{28.67 9} {29.33 10}
  \arr{32.67 9} {33.33 10}
  \arr{22.67 10} {23.33 9}
  \arr{26.67 10} {27.33 9}
  \arr{30.67 10} {31.33 9}
  \arr{34.67 10} {35.23 9.15}
  \arr{20.67 13} {21.33 12}
  \arr{24.67 13} {25.33 12}
  \arr{28.67 13} {29.33 12}
  \arr{32.67 13} {33.33 12}
  \arr{36.67 13} {37.33 12}
  \arr{18.67 12} {19.33 13}
  \arr{22.67 12} {23.33 13}
  \arr{26.67 12} {27.33 13}
  \arr{30.67 12} {31.33 13}
  \arr{34.67 12} {35.33 13}
  \arr{16.67 15} {17.33 16}
  \arr{20.67 15} {21.33 16}
  \arr{24.67 15} {25.33 16}
  \arr{28.67 15} {29.33 16}
  \arr{32.67 15} {33.33 16}
  \arr{36.67 15} {37.33 16}
  \arr{18.67 16} {19.33 15}
  \arr{22.67 16} {23.33 15}
  \arr{26.67 16} {27.33 15}
  \arr{30.67 16} {31.33 15}
  \arr{34.67 16} {35.33 15}
  \arr{38.67 16} {39.33 15}
  \arr{16.5 19.25} {17.33 18}
  \arr{20.5 19.25} {21.33 18}
  \arr{24.5 19.25} {25.33 18}
  \arr{28.5 19.25} {29.33 18}
  \arr{32.5 19.25} {33.33 18}
  \arr{36.5 19.25} {37.33 18}
  \arr{40.5 19.25} {41.33 18}
  \arr{14.67 18} {15.5 19.25}
  \arr{18.67 18} {19.5 19.25}
  \arr{22.67 18} {23.5 19.25}
  \arr{26.67 18} {27.5 19.25}
  \arr{30.67 18} {31.5 19.25}
  \arr{34.67 18} {35.5 19.25}
  \arr{38.67 18} {39.5 19.25}
  \arr{24.5 20.75} {25.5 22.25}
  \arr{26.5 22.25} {27.5 20.75}
  \setdots<2pt>
  \plot 20.5 -.15  23.5 -.15 /
  \plot 26.5 -4.15  29.5 -4.15 /
  \plot 30.5 -4.15  33.5 -4.15 /
  \plot 14   19.85  15.5 19.85 /
  \plot 16.5 19.85  19.5 19.85 /
  \plot 20.5 19.85  23.5 19.85 /
  \plot 28.5 19.85  31.5 19.85 /
  \plot 32.5 19.85  35.5 19.85 /
  \plot 36.5 19.85  39.5 19.85 /
  \plot 40.5 19.85  42   19.85 /
  \setshadegrid span <1mm>
  \vshade 23   -.5 0 <,z,,>
          23.5 -1.5 1 <z,z,,>
          24.25 -1.5 3 <z,z,,>
          27    4 8 <z,z,,>
          28    6 9.5 <z,z,,>
          30    5.5 12.5 <z,z,,>
          30.8  6.8 13.7 <z,z,,>
          31    11  14 <z,z,,>
          35.67 18 21 <z,z,,>
          36    18.5 21 <z,z,,>
          37    20 20.5 /
  \multiput{$*$} at 28 -3  32 -3  18 19  22 19  26 20  30 19  34 19  38 19  /
  \put{$(1)$} at 32 17
  \put{$(2)$} at 30 14
  \put{$(3)$} at 28 11
  \put{$(4)$} at 28 4
  \put{$(5)$} at 26 0
  \endpicture$$
\end{figure}

In Example~\ref{bounded-I} we have seen that the sequences labelled 
by a $*$ are Auslander-Reiten sequences.
It is straightforward to verify that also the other meshes in the picture correspond to
nonsplit short exact sequences.  We show that each sequence is almost split.
From the sequences labelled by a $*$ we obtain that the first term in a sequence 
given by one of the meshes $(1)$ or $(5)$ is the translate of the last term. 
Since every nonautomorphism of
$C_1=\vcenter{\beginpicture 
  \setcoordinatesystem units <.3cm,.3cm>
    \multiput{$\smallsq5$} at 0 0  0 .5 /
    \multiput{} at 0 -.5  .5 1 /
    \put{\bulletsymb} at .25 0
    \endpicture}$
  factors through
$\vcenter{\beginpicture 
    \setcoordinatesystem units <.3cm,.3cm>
    \multiput{$\smallsq5$} at 0 0 /
    \multiput{} at 0 -.5  .5 1 /
    \put{\bulletsymb} at .25 0
    \endpicture}\to
  \vcenter{\beginpicture 
    \setcoordinatesystem units <.3cm,.3cm>
    \multiput{$\smallsq5$} at 0 0  0 .5 /
    \multiput{} at 0 -.5  .5 1 /
    \put{\bulletsymb} at .25 0
    \endpicture}$,
  the sequence labelled $(1)$ is an Auslander-Reiten sequence, by Proposition~\ref{ar-test},
  $(2.\Rightarrow1.)$.
  As a consequence, the sequences ending at $C_1$, $\tau C_1$, $\tau^2 C_1$, $\tau^3 C_1$ 
  and $\tau^{-1}C_1$
  are all Auslander-Reiten sequences
  with the additional property that their middle term has 
  two indecomposable summands which are as pictured.
  Similarly, every nonautomorphism of 
  $A_5=\vcenter{\beginpicture 
  \setcoordinatesystem units <.3cm,.3cm>
    \multiput{$\smallsq5$} at 0 0  0 .5  0 1 /
    \multiput{} at 0 -.5  .5 1 /
    \put{\bulletsymb} at .25 .5
    \endpicture}$
  factors over the inclusion
  $\vcenter{\beginpicture 
  \setcoordinatesystem units <.3cm,.3cm>
    \multiput{$\smallsq5$} at 0 0  0 .5  0 1 /
    \multiput{} at 0 -.5  .5 1 /
    \put{\bulletsymb} at .25 .5
    \endpicture}\to
  \vcenter{\beginpicture 
  \setcoordinatesystem units <.3cm,.3cm>
    \multiput{$\smallsq5$} at 0 0  0 .5  0 1 /
    \multiput{} at 0 -.5  .5 1 /
    \put{\bulletsymb} at .25 1
    \endpicture}$,
  hence it follows from Proposition~\ref{ar-test} that the sequence labelled $(5)$
  is an Auslander-Reiten sequence with middle term a direct sum of two indecomposables,
  as in the picture.  This holds also true for the sequence starting at $\tau^{-1}A_5$. 
  In order to verify that the sequences labelled $(2)$, $(3)$, $(4)$ are almost split,
  consider the modules 
  $$C_2=  \vcenter{\beginpicture 
  \setcoordinatesystem units <.3cm,.3cm>
    \multiput{$\smallsq5$} at 0 0  0 .5  0 1  0 1.5  0 2  0 2.5  .5 .5  .5 1 /
    \multiput{} at 0 -.5  .5 1 /
    \multiput{\bulletsymb} at .25 .5  .75 .5 /
    \plot  .25 .5  .75 .5 /
    \endpicture}\;; \quad
  C_3=  \vcenter{\beginpicture 
  \setcoordinatesystem units <.3cm,.3cm>
    \multiput{$\smallsq5$} at 0 0  0 .5  0 1  0 1.5  0 2  .5 .5  .5 1 /
    \multiput{} at 0 -.5  .5 1 /
    \multiput{\bulletsymb} at .25 .5  .75 .5 /
    \plot  .25 .5  .75 .5 /
    \endpicture}\;; \quad
  A_4=  \vcenter{\beginpicture 
  \setcoordinatesystem units <.3cm,.3cm>
    \multiput{$\smallsq5$} at 0 0  0 .5  0 1  0 1.5  .5 .5  .5 1 /
    \multiput{} at 0 -.5  .5 1 /
    \multiput{\bulletsymb} at .25 .5  .75 .5 /
    \plot  .25 .5  .75 .5 /
    \endpicture}\;.
  $$
  Let $X$ be one of these modules, then $X=(U\subset V)$ where $V=V_1\oplus V_2$
  with $V_1=\Lambda/\rad^m\!\Lambda$ with $m=4,5,6$ and $V_2=\Lambda/\rad^2\!\Lambda$.
  Identifying $\End_{\mathcal S}X$ as a subset of $\End_\Lambda V$, it is easy to see
  that the radical of the endomorphism ring of the object $X$ is given as follows.
  $$\rad\End_{\mathcal S}X\quad=\quad\rad\End_\Lambda V\quad=\quad
  \left({\rad\Lambda/\rad^m\!\Lambda \atop \Lambda/\rad^2\!\Lambda}\;
    {\Lambda/\rad^2\!\Lambda \atop \rad\Lambda/\rad^2\!\Lambda}\right)$$
\def\arqthreesix{%
\beginpicture\setcoordinatesystem units <4.56mm,4.5mm>
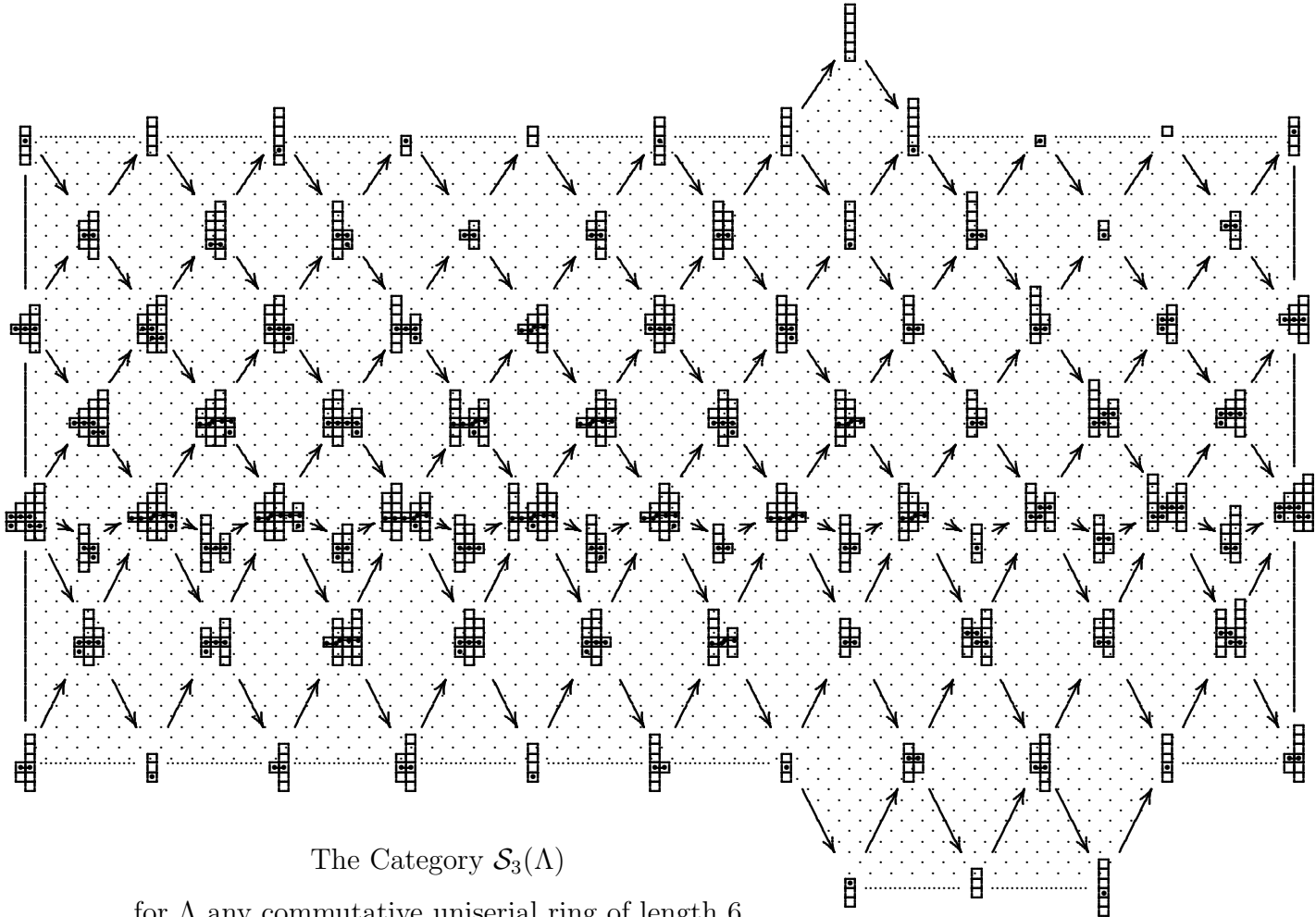
\put{The Category $\mathcal S_3(\Lambda)$} at 13 -3.3
\put{for $\Lambda$ any commutative uniserial ring of length 6} at 13 -5
\put{} at 42 0
\put{} at -3 0
\multiput{$\smallsq3$} at  25.85 -4.6  25.85 -4.3  25.85 -4 /
\multiput{\bulletsymb} at 26 -4 /

\multiput{$\smallsq3$} at 29.85 -4.3  29.85 -4  29.85 -3.7 /

\multiput{$\smallsq3$} at 33.85 -4.9  33.85 -4.6  33.85 -4.3  33.85 -4
        33.85 -3.7  33.85 -3.4 /
\multiput{\bulletsymb} at 34 -4.3 /
  \multiput{$\smallsq3$} at 0 0  0 .3  0 .6  0 -.3  0 -.6  0 -.9  -.3 -.3  
  -.3 -.6 /
  \multiput{\bulletsymb} at -.15 -.3  .15 -.3 /
  \plot -.15 -.3  .15 -.3 /
  
  \multiput{$\smallsq3$} at 3.85 0  3.85 -.3  3.85 -.6 /
  \multiput{\bulletsymb} at 4 -.6 /
  
  \multiput{$\smallsq3$} at 7.7 -.3  8 -.9  8 -.6  8 -.3  8 0  8 .3 /
  \multiput{\bulletsymb} at 7.85 -.3  8.15 -.3 /
  \plot  7.85 -.3  8.15 -.3 /
  
  \multiput{$\smallsq3$} at 11.7 -.6  11.7 -.3  11.7 0  12 -.9  12 -.6  12 -.3
  12 0  12 .3  12 .6 /
  \multiput{\large\bf.} at 11.85 -.3  12.15 -.3 /
  \plot 11.85 -.3  12.15 -.3 /
  
  \multiput{$\smallsq3$} at 15.85 -.6  15.85 -.3  15.85 0  15.85 .3 /
  \multiput{\bulletsymb} at 16 -.6 /
  
  \multiput{$\smallsq3$} at 19.7 -.9   19.7 -.6  19.7 -.3  19.7 0  19.7 .3
  19.7 .6  20 -.3 /
  \multiput{\bulletsymb} at 19.85 -.3  20.15 -.3 /
  \plot 19.85 -.3  20.15 -.3 /
  
  \multiput{$\smallsq3$} at 23.85 -.6  23.85 -.3  23.85 0  /
  \multiput{\bulletsymb} at 24 -.3 /
  
  \multiput{$\smallsq3$} at 27.7 -.6  27.7 -.3  27.7 0  27.7 .3 28 -.3  28 0 /
  \multiput{\bulletsymb} at 27.85 0  28.15 0 /
  \plot 27.85 0  28.15 0 /
  
  \multiput{$\smallsq3$} at 31.7 -.3  31.7 0  31.7 .3  32 -.9  32 -.6  32 -.3 
  32 0  32 .3  32 .6 /
  \multiput{\bulletsymb} at 31.85 -.3  32.15 -.3 /
  \plot 31.85 -.3  32.15 -.3 /
  
  \multiput{$\smallsq3$} at 35.85 -.6  35.85 -.3  35.85 0  35.85 .3  35.85 .6 /
  \multiput{\bulletsymb} at 36 -.3 /
  
  \multiput{$\smallsq3$} at 39.7 -.3  39.7 0  40 -.6  40 -.3  40 0  40 .3 
  40 .6  40 .9 /
  \multiput{\bulletsymb} at 39.85 0  40.15 0 /
  \plot 39.85 0  40.15 0 /
  \multiput{$\smallsq3$} at 1.55 3.4  1.55 3.7  1.85 3.1  1.85 3.4  1.85 3.7
  1.85 4  1.85 4.3  1.85 4.6  2.15 3.4  2.15 3.7  2.15 4 /
  \multiput{\bulletsymb} at 1.7 3.4  1.7 3.7  2 3.7  2.3 3.7 /
  \plot 1.7 3.7  2.3 3.7 /
  
  \multiput{$\smallsq3$} at 5.55 3.4  5.55 3.7  5.55 4  5.85 3.7  6.15 3.1 
  6.15 3.4  6.15 3.7  6.15 4  6.15 4.3 /
  \multiput{\bulletsymb} at 5.7 3.4  5.7 3.7  6 3.7  6.3 3.7 /
  \plot 5.7 3.7  6.3 3.7 /
  
  \multiput{$\smallsq3$} at  9.4 3.7  9.7 3.1  9.7 3.4  9.7 3.7  9.7 4
  9.7 4.3  10 3.4  10 3.7  10 4  10.3 3.1  10.3 3.4  10.3 3.7 
  10.3 4  10.3 4.3  10.3 4.6 /
  \multiput{\bulletsymb} at 9.5 3.65  9.8 3.65  9.9 3.75  10.2 3.75  
  10.5 3.75 / 
  \plot 9.5 3.65  9.8 3.65 /
  \plot 9.9 3.75  10.5 3.75 /
  
  \multiput{$\smallsq3$} at 13.55 3.4  13.55 3.7  13.55 4  13.85 3.1  
  13.85 3.4  13.85 3.7  13.85 4  13.85 4.3  13.85 4.6
  14.15 3.4  14.15 3.7  14.15 4  14.15 4.3 /
  \multiput{\bulletsymb} at 13.7 3.4  13.7 3.7  14 3.7  14.3 3.7 /
  \plot 13.7 3.7  14.3 3.7 /
  
  \multiput{$\smallsq3$} at 17.55 3.4  17.55 3.7  17.55 4  17.55 4.3 
  17.85 3.1  17.85 3.4  17.85 3.7  17.85 4  17.85 4.3  17.85 4.6
  18.15 3.7 /
  \multiput{\bulletsymb} at 17.7 3.4  17.7 3.7  18 3.7  18.3 3.7 /
  \plot  17.7 3.7  18.3 3.7 /
  
  \multiput{$\smallsq3$} at 21.55 3.1  21.55 3.4  21.55 3.7  21.55 4  21.55 4.3
  21.55 4.6  21.85 3.7  22.15 3.4  22.15 3.7  22.15 4 /
  \multiput{\bulletsymb} at 21.65 3.65  21.95 3.65  22.05 3.75  22.35 3.75 /
  \plot 21.65 3.65  21.95 3.65 /
  \plot  22.05 3.75  22.35 3.75 /
  
  \multiput{$\smallsq3$} at 25.7 3.4  25.7 3.7  25.7 4  25.7 4.3  26 3.7  26 4 /
  \multiput{\bulletsymb} at 25.85 3.7  26.15 3.7 /
  \plot  25.85 3.7  26.15 3.7 /
  
  \multiput{$\smallsq3$} at 29.55 3.4  29.55 3.7  29.55 4  29.55 4.3  29.85 3.7
  29.85 4  30.15 3.1  30.15 3.4  30.15 3.7  30.15 4  30.15 4.3
  30.15 4.6 /
  \multiput{\bulletsymb} at 29.7 4  30 4  30 3.7  30.3 3.7 /
  \plot 29.7 4  30 4 /
  \plot  30 3.7  30.3 3.7 /
  
  \multiput{$\smallsq3$} at 33.7 3.7  33.7 4  33.7 4.3  34 3.4  34 3.7  34 4
  34 4.3  34 4.6 /
  \multiput{\bulletsymb} at 33.85 3.7  34.15 3.7 /
  \plot 33.85 3.7  34.15 3.7 /
  
  \multiput{$\smallsq3$} at 37.55 3.4  37.55 3.7  37.55 4  37.55 4.3  37.55 4.6
  37.85 3.7  37.85 4  38.15 3.4  38.15 3.7  38.15 4  38.15 4.3
  38.15 4.6  38.15 4.9 /
  \multiput{\bulletsymb} at 37.7 4  38 4  38 3.7  38.3 3.7 /
  \plot 37.7 4  38 4 /
  \plot  38 3.7  38.3 3.7 /
  
  \multiput{$\smallsq3$} at -.6 7.4  -.6 7.7  -.3 7.4  -.3 7.7  -.3 8
  0 7.1  0 7.4  0 7.7  0 8  0 8.3  .3 7.1  .3 7.4  .3 7.7  .3 8
  .3 8.3  .3 8.6 /
  \multiput{\bulletsymb} at  -.45 7.4  -.45 7.7  -.15 7.7  .15 7.7  .15 7.4
  .45 7.4 /
  \plot -.45 7.7  .15 7.7 /
  \plot .15 7.4   .45 7.4 /
  
  \multiput{$\smallsq3$} at 3.25 7.7  3.55 7.4  3.55 7.7  3.55 8  
  3.85 7.1  3.85 7.4  3.85 7.7  3.85 8  3.85 8.3
  4.15 7.1  4.15 7.4  4.15 7.7  4.15 8  4.15 8.3  4.15 8.6
  4.45 7.4  4.45 7.7 /
  \multiput{\bulletsymb} at 3.35 7.65  3.65 7.65  3.95 7.65  
  4.05 7.75  4.35 7.75  4.65 7.75  4.6 7.4 /
  \plot 3.35 7.65  3.95 7.65 /
  \plot 4.05 7.75  4.65 7.75 /
  
  \multiput{$\smallsq3$} at 7.25 7.4  7.25 7.7  7.25 8  7.55 7.1  7.55 7.4
  7.55 7.7  7.55 8  7.55 8.3  7.85 7.1  7.85 7.4  7.85 7.7
  7.85 8  7.85 8.3  7.85 8.6  8.15 7.7  8.45 7.4  8.45 7.7
  8.45 8 /
  \multiput{\bulletsymb} at 7.35 7.65  7.65 7.65  7.75 7.75  8.05 7.75
  8.35 7.75  8.65 7.75  8.6 7.4 /
  \plot 7.35 7.65  7.65 7.65 /
  \plot 7.75 7.75  8.65 7.75 /
  
  \multiput{$\smallsq3$} at 11.25 7.4  11.25 7.7  11.25 8  11.25 8.3
  11.55 7.1  11.55 7.4  11.55 7.7  11.55 8  11.55 8.3  11.55 8.6
  11.85 7.7  12.15 7.4  12.15 7.7  12.15 8  12.45 7.1  12.45 7.4
  12.45 7.7  12.45 8  12.45 8.3 / 
  \multiput{\bulletsymb} at 11.35 7.65  11.65 7.65  11.95 7.65  12.25 7.65
  12.35 7.75  12.65 7.75  12.3  7.4 /
  \plot 11.35 7.65  12.25 7.65 /
  \plot 12.35 7.75  12.65 7.75 /
  
  \multiput{$\smallsq3$} at 15.25 7.1  15.25 7.4  15.25 7.7  15.25 8  15.25 8.3
  15.25 8.6  15.55 7.7  15.85 7.4  15.85 7.7  15.85 8  
  16.15 7.1  16.15 7.4  16.15 7.7  16.15 8  16.15 8.3  16.15 8.6
  16.45 7.4  16.45 7.7  16.45 8  16.45 8.3  /
  \multiput{\bulletsymb} at 15.35 7.65  15.65 7.65  15.95 7.65 
  16.05 7.75  16.35 7.75  16.65 7.75  16 7.4 /
  \plot 15.35 7.65  15.95 7.65 /
  \plot 16.05 7.75  16.65 7.75 /
  
  \multiput{$\smallsq3$} at 19.4 7.7  19.7 7.4  19.7 7.7  19.7 8  20 7.1
  20 7.4  20 7.7  20 8  20 8.3  20 8.6  20.3 7.4  20.3 7.7
  20.3 8  20.3 8.3 /
  \multiput{\bulletsymb} at 19.5 7.65  19.8 7.65  19.9 7.75  20.2 7.75  20.5 7.75
  20.45 7.4 /
  \plot 19.5 7.65  19.8 7.65 /
  \plot 19.9 7.75  20.5 7.75 /
  
  \multiput{$\smallsq3$} at 23.4 7.7  23.4 8  23.7 7.1  23.7 7.4  23.7 7.7
  23.7 8  23.7 8.3  23.7 8.6  24 7.4  24 7.7  24 8  24 8.3 
  24.3 7.7 /
  \multiput{\bulletsymb} at 23.5 7.65  23.8 7.65  24.1 7.65  24.2 7.75  
  24.5 7.75 /
  \plot 23.5 7.65  24.1 7.65 /
  \plot 24.2 7.75  24.5 7.75 /
  
  \multiput{$\smallsq3$} at 27.55 7.1  27.55 7.4  27.55 7.7  27.55 8  27.55 8.3
  27.55 8.6  27.85 7.4  27.85 7.7  27.85 8  27.85 8.3  28.15 7.7  28.15 8         /
  \multiput{\bulletsymb} at 27.65 7.65  27.95 7.65  28.05 7.75  28.35 7.75 /
  \plot 27.65 7.65  27.95 7.65 /
  \plot   28.05 7.75  28.35 7.75 /
  
  \multiput{$\smallsq3$} at 31.55 7.4  31.55 7.7  31.55 8  31.55 8.3  31.55 8.6
  31.85 7.7  31.85 8  32.15 7.4  32.15 7.7  32.15 8  32.15 8.3 
  /
  \multiput{\bulletsymb} at 31.7 7.7  32 7.7  32 8 32.3 8 /
  \plot 31.7 7.7  32 7.7 /
  \plot  32 8 32.3 8 /
  
  \multiput{$\smallsq3$} at   35.4 7.4  35.4 7.7  35.4 8  35.4 8.3  35.4 8.6  
  35.4 8.9  35.7 7.7  35.7 8  36 7.7  36 8  36 8.3
  36.3 7.4  36.3 7.7  36.3 8  36.3 8.3  36.3 8.6 /
  \multiput{\bulletsymb} at 35.55 7.7  35.85 7.7  35.85 8  36.15 8  36.45 8 /
  \plot 35.55 7.7  35.85 7.7  /
  \plot 35.85 8  36.45 8 /
  
  \multiput{$\smallsq3$} at  39.4 7.7  39.4 8  39.7 7.7  39.7 8  39.7 8.3 
  40 7.4  40 7.7  40 8  40 8.3  40 8.6  
  40.3 7.4  40.3 7.7  40.3 8  40.3 8.3  40.3 8.6  40.3 8.9 /
  \multiput{\bulletsymb} at 39.55 7.7  39.55 8  39.85 8  40.15 8  
  40.15 7.7  40.45 7.7 /
  \plot 39.55 8  40.15 8 /
  \plot 40.15 7.7  40.45 7.7 /
  
  \multiput{$\smallsq3$} at 1.7 6.1  1.7 6.4  1.7 6.7  1.7 7  1.7 7.3 
  2 6.4  2 6.7 /
  \multiput{\bulletsymb} at 1.85 6.7  2.15 6.7  2.15 6.4 /
  \plot 1.85 6.7  2.15 6.7 /
  
  \multiput{$\smallsq3$} at 5.55 6.1  5.55 6.4  5.55 6.7  5.55 7  5.55 7.3
  5.55 7.6  5.85 6.7  6.15 6.4  6.15 6.7  6.15 7 /
  \multiput{\bulletsymb} at 5.7 6.7  6 6.7  6.3 6.7 /
  \plot 5.7 6.7  6.3 6.7 /
  
  \multiput{$\smallsq3$} at 9.7 6.4  9.7 6.7  9.7 7  10 6.1  10 6.4  10 6.7
  10 7 10 7.3 /
  \multiput{\bulletsymb} at 9.85 6.7  10.15 6.7  9.85 6.4 /
  \plot 9.85 6.7  10.15 6.7 /
  
  \multiput{$\smallsq3$} at 13.55 6.1  13.55 6.4  13.55 6.7  13.55 7  13.55 7.3 
  13.55 7.6  13.85 6.4  13.85 6.7  13.85 7  13.85 7.3  14.15 6.7
  /
  \multiput{\bulletsymb} at 13.7 6.7  14 6.7  14.3 6.7 /
  \plot 13.7 6.7  14.3 6.7 /
  
  \multiput{$\smallsq3$} at 17.7 6.1  17.7 6.4  17.7 6.7  17.7 7  17.7 7.3 
  17.7 7.6  18 6.4  18 6.7  18 7 /
  \multiput{\bulletsymb} at 17.85 6.7  18.15 6.7  18.15 6.4 /
  \plot 17.85 6.7  18.15 6.7 /
  
  \multiput{$\smallsq3$} at 21.7 6.4  21.7 6.7  21.7 7  21.7 7.3  22 6.7 /
  \multiput{\bulletsymb} at 21.85 6.7  22.15 6.7 /
  \plot 21.85 6.7  22.15 6.7 /
  
  \multiput{$\smallsq3$} at 25.7 6.1  25.7 6.4  25.7 6.7  25.7 7  25.7 7.3 
  25.7 7.6  26 6.7  26 7 /
  \multiput{\bulletsymb} at 25.85 6.7  26.15 6.7 /
  \plot  25.85 6.7  26.15 6.7 /
  
  \multiput{$\smallsq3$} at 29.85 6.4  29.85 6.7  29.85 7  29.85 7.3 /
  \multiput{\bulletsymb} at 30 6.7 /
  
  \multiput{$\smallsq3$} at 33.7 6.4  33.7 6.7  33.7 7  33.7 7.3  33.7 7.6 
  34 6.7  34 7 /
  \multiput{\bulletsymb} at 33.85 7  34.15 7 /
  \plot 33.85 7  34.15 7 /
  
  \multiput{$\smallsq3$} at 37.7 6.7  37.7 7  37.7 7.3  38 6.4  38 6.7
  38 7  38 7.3  38 7.6  38 7.9 /
  \multiput{\bulletsymb} at 37.85 6.7  38.15 6.7 /
  \plot 37.85 6.7  38.15 6.7 /
  
  \multiput{$\smallsq3$} at 1.4 10.7  1.7 10.4  1.7 10.7  1.7 11  2 10.1
  2 10.4  2 10.7  2 11  2 11.3  2.3 10.1  2.3 10.4  2.3 10.7
  2.3 11  2.3 11.3  2.3 11.6 /
  \multiput{\bulletsymb} at 1.55 10.7  1.85 10.7  2.15 10.7  2.15 10.4  2.45 10.4
  /
  \plot 1.55 10.7  2.15 10.7 /
  \plot 2.15 10.4  2.45 10.4 /
  
  \multiput{$\smallsq3$} at 5.4 10.4  5.4 10.7  5.4 11  5.7 10.1  5.7 10.4
  5.7 10.7  5.7 11  5.7 11.3  6 10.1  6 10.4  6 10.7  6 11
  6 11.3  6 11.6  6.3 10.4  6.3 10.7 /
  \multiput{\bulletsymb} at 5.5 10.65  5.8 10.65  5.9 10.75  6.2 10.75
  6.5 10.75  6.45 10.4 /
  \plot 5.5 10.65  5.8 10.65 /
  \plot 5.9 10.75  6.5 10.75 /
  
  \multiput{$\smallsq3$} at 9.4 10.4  9.4 10.7  9.4 11  9.4 11.3  9.7 10.1
  9.7 10.4  9.7 10.7  9.7 11 9.7 11.3  9.7 11.6  10 10.7
  10.3 10.4  10.3 10.7  10.3 11 /
  \multiput{\bulletsymb} at 9.55 10.7  9.85 10.7  10.15 10.7  10.45 10.7
  10.45 10.4 /
  \plot 9.55 10.7  10.45 10.7 /
  
  \multiput{$\smallsq3$} at 13.4 10.1  13.4 10.4  13.4 10.7  13.4 11  13.4 11.3
  13.4 11.6  13.7 10.7  14 10.4  14 10.7  14 11  14.3  10.1
  14.3 10.4  14.3 10.7  14.3 11  14.3 11.3 /
  \multiput{\bulletsymb} at 13.5 10.65  13.8 10.65  14.1 10.65  14.2 10.75
  14.5 10.75  14.15 10.4 /
  \plot 13.5 10.65  14.1 10.65 /
  \plot 14.2 10.75  14.5 10.75 /
  
  \multiput{$\smallsq3$} at 17.4 10.7  17.7 10.4  17.7 10.7  17.7 11  
  18 10.1  18 10.4  18 10.7  18 11  18 11.3  18 11.6  
  18.3 10.4  18.3 10.7  18.3 11  18.3 11.3 /
  \multiput{\bulletsymb} at 17.5 10.65  17.8 10.65  17.9 10.75  18.2 10.75
  18.5 10.75 /
  \plot 17.5 10.65  17.8 10.65 /
  \plot 17.9 10.75  18.5 10.75 /
  
  \multiput{$\smallsq3$} at 21.55 10.7  21.55 11  21.85 10.1  21.85 10.4
  21.85 10.7  21.85 11  21.85 11.3  21.85 11.6  
  22.15 10.4  22.15 10.7  22.15 11  22.15 11.3 /
  \multiput{\bulletsymb} at 21.7 10.7  22 10.7  22.3 10.7  22.3 10.4 /
  \plot 21.7 10.7  22.3 10.7 /
  
  \multiput{$\smallsq3$} at 25.55 10.1  25.55 10.4  25.55 10.7  25.55 11 
  25.55 11.3  25.55 11.6  25.85 10.4  25.85 10.7  25.85 11
  25.85 11.3  26.15 10.7 / 
  \multiput{\bulletsymb} at 25.65 10.65  25.95 10.65  26.05 10.75  26.35 10.75 /
  \plot 26.05 10.75  26.35 10.75 /
  \plot 25.65 10.65  25.95 10.65 /
  
  \multiput{$\smallsq3$} at 29.7 10.4  29.7 10.7  29.7 11 29.7 11.3  29.7 11.6 
  30 10.7  30 11 /
  \multiput{\bulletsymb} at 29.85 10.7  30.15 10.7 /
  \plot 29.85 10.7  30.15 10.7 /
  
  \multiput{$\smallsq3$} at 33.55 10.4  33.55 10.7  33.55 11  33.55 11.3  
  33.55 11.6  33.55 11.9  33.85 10.7  33.85 11  34.15 10.4  34.15 10.7
  34.15 11  34.15 11.3 /
  \multiput{\bulletsymb} at 33.7 10.7  34 10.7  34 11  34.3 11 /
  \plot 34 11  34.3 11 /
  \plot  33.7 10.7  34 10.7 /
  
  \multiput{$\smallsq3$} at 37.55 10.7  37.55 11  37.85 10.7  37.85 11  
  37.85 11.3  38.15 10.4  38.15 10.7  38.15 11  38.15 11.3  38.15 11.6 /
  \multiput{\bulletsymb} at 37.7 10.7  37.7 11  38 11  38.3 11 /
  \plot 37.7 11  38.3 11 /
  
  \multiput{$\smallsq3$} at -.45 13.7  -.15 13.4  -.15 13.7
  -.15 14  .15 13.1  .15 13.4  .15 13.7  .15 14  .15 14.3 /
  \multiput{\bulletsymb} at -.3 13.7  0 13.7  .3 13.7 /
  \plot -.3 13.7  .3 13.7 /
  
  \multiput{$\smallsq3$} at 3.55 13.4  3.55 13.7  3.55 14  3.85 13.1  
  3.85 13.4  3.85 13.7  3.85 14  3.85 14.3  4.15 13.1  4.15 13.4 
  4.15 13.7  4.15 14  4.15 14.3  4.15 14.6 /
  \multiput{\bulletsymb} at 3.7 13.7  4 13.7  4 13.4  4.3 13.4 /
  \plot 4 13.4  4.3 13.4 /
  \plot 3.7 13.7  4 13.7 /
  
  \multiput{$\smallsq3$} at 7.55 13.4  7.55 13.7  7.55 14  7.55 14.3  7.85 13.1
  7.85 13.4  7.85 13.7  7.85 14  7.85 14.3  7.85 14.6 
  8.15 13.4  8.15 13.7 /
  \multiput{\bulletsymb} at 7.7 13.7  8 13.7  8.3 13.7  8.3 13.4 /
  \plot  7.7 13.7  8.3 13.7 /
  
  \multiput{$\smallsq3$} at 11.55 13.1  11.55 13.4  11.55 13.7  11.55 14 
  11.55 14.3  11.55 14.6  11.85 13.7  12.15 13.4  12.15 13.7  12.15 14 /
  \multiput{\bulletsymb} at 11.7 13.7  12 13.7  12.3 13.7  12.3 13.4 /
  \plot 11.7 13.7  12.3 13.7 /
  
  \multiput{$\smallsq3$} at 15.55 13.7  15.85 13.4  15.85 13.7  15.85 14
  16.15 13.1  16.15 13.4  16.15 13.7  16.15 14  16.15 14.3 /
  \multiput{\bulletsymb} at 15.65 13.65  15.95 13.65  16.05 13.75  16.35 13.75 /
  \plot  16.05 13.75  16.35 13.75 /
  \plot  15.65 13.65  15.95 13.65 /
  
  \multiput{$\smallsq3$} at 19.55 13.7  19.55 14  19.85 13.1  19.85 13.4  
  19.85 13.7  19.85 14  19.85 14.3  19.85 14.6  20.15 13.4  20.15 13.7 
  20.15 14  20.15 14.3 /
  \multiput{\bulletsymb} at 19.7 13.7  20 13.7  20.3 13.7 /
  \plot 19.7 13.7  20.3 13.7 /
  
  \multiput{$\smallsq3$} at 23.7 13.1  23.7 13.4  23.7 13.7  23.7 14 
  23.7 14.3  23.7 14.6  24 13.4  24 13.7  24 14  24 14.3 /
  \multiput{\bulletsymb} at 23.85 13.7  24.15 13.7  24.15 13.4 /
  \plot 23.85 13.7  24.15 13.7 /
  
  \multiput{$\smallsq3$} at 27.7 13.4  27.7 13.7  27.7 14  27.7 14.3  27.7 14.6
  28 13.7 /
  \multiput{\bulletsymb} at 27.85 13.7  28.15 13.7 /
  \plot 27.85 13.7  28.15 13.7 /
  
  \multiput{$\smallsq3$} at 31.7 13.4  31.7 13.7  31.7 14  31.7 14.3  31.7 14.6
  31.7 14.9  32 13.7  32 14 /
  \multiput{\bulletsymb} at 31.85 13.7  32.15 13.7 /
  \plot 31.85 13.7  32.15 13.7 /
  
  \multiput{$\smallsq3$} at 35.7 13.7  35.7 14  36 13.4  36 13.7  36 14  36 14.3 
  /
  \multiput{\bulletsymb} at 35.85 13.7  35.85 14  36.15 14 /
  \plot 35.85 14  36.15 14 /
  
  \multiput{$\smallsq3$} at 39.55 14  39.85 13.7  39.85 14  
  39.85 14.3  40.15 13.4  40.15 13.7  40.15 14  40.15 14.3  40.15 14.6 /
  \multiput{\bulletsymb} at 39.7 14  40 14  40.3 14 /
  \plot 39.7 14  40.3 14 /
  
  \multiput{$\smallsq3$} at 1.7 16.4  1.7 16.7  1.7 17
  2 16.1  2 16.4  2 16.7  2 17  2 17.3 /
  \multiput{\bulletsymb} at 1.85 16.7  2.15 16.7 /
  \plot  1.85 16.7  2.15 16.7 /
  
  \multiput{$\smallsq3$} at 5.7 16.4  5.7 16.7  5.7 17  5.7 17.3  6 16.1
  6 16.4  6 16.7  6 17  6 17.3  6 17.6 /
  \multiput{\bulletsymb} at 5.85 16.4  6.15 16.4 /
  \plot 5.85 16.4  6.15 16.4 /
  
  \multiput{$\smallsq3$} at 9.7 16.1  9.7 16.4  9.7 16.7  9.7 17  9.7 17.3
  9.7 17.6  10 16.4  10 16.7 /
  \multiput{\bulletsymb} at 9.85 16.7  10.15 16.7  10.15 16.4 /
  \plot 9.85 16.7  10.15 16.7 /
  
  \multiput{$\smallsq3$} at 13.7 16.7  14 16.4  14 16.7  14 17 /
  \multiput{\bulletsymb} at 13.85 16.7  14.15 16.7 /
  \plot 13.85 16.7  14.15 16.7 /
  
  \multiput{$\smallsq3$} at 17.7 16.7  17.7 17  18 16.1  18 16.4  18 16.7
  18 17  18 17.3 /
  \multiput{\bulletsymb} at 17.85 16.7  18.15 16.7 /
  \plot 17.85 16.7  18.15 16.7 /
  
  \multiput{$\smallsq3$} at 21.7 16.1  21.7 16.4  21.7 16.7  21.7 17  21.7 17.3
  21.7 17.6  22 16.4  22 16.7  22 17  22 17.3 /
  \multiput{\bulletsymb} at 21.85 16.7  22.15 16.7 /
  \plot 21.85 16.7  22.15 16.7 /
  
  \multiput{$\smallsq3$} at 25.85 16.4  25.85 16.7  25.85 17  25.85 17.3
  25.85 17.6 /
  \multiput{\bulletsymb} at 26 16.4 /
  
  \multiput{$\smallsq3$} at 29.7 16.4  29.7 16.7  29.7 17  29.7 17.3  29.7 17.6
  29.7 17.9  30 16.7 /
  \multiput{\bulletsymb} at 29.85 16.7  30.15 16.7 /
  \plot  29.85 16.7  30.15 16.7 /
  
  \multiput{$\smallsq3$} at 33.85 16.7  33.85 17 /
  \multiput{\bulletsymb} at 34 16.7 /
  
  \multiput{$\smallsq3$} at 37.7 17  38 16.4  38 16.7  38 17  38 17.3 /
  \multiput{\bulletsymb} at 37.85 17  38.15 17 /
  \plot 37.85 17  38.15 17 /
  
  \multiput{$\smallsq3$} at -.15 19.1  -.15 19.4  -.15 19.7  -.15 20 /
  \multiput{\bulletsymb} at 0 19.7 /
  
  \multiput{$\smallsq3$} at 3.85 19.4  3.85 19.7  3.85 20  3.85 20.3 /
  
  \multiput{$\smallsq3$} at 7.85 19.1  7.85 19.4  7.85 19.7  7.85 20 
  7.85 20.3  7.85 20.6 /
  \multiput{\bulletsymb} at 8 19.4 /
  
  \multiput{$\smallsq3$} at 11.85 19.4  11.85 19.7 /
  \multiput{\bulletsymb} at 12 19.7 /
  
  \multiput{$\smallsq3$} at 15.85 19.7  15.85 20 /
  
  \multiput{$\smallsq3$} at 19.85 19.1  19.85 19.4  19.85 19.7  19.85 20
  19.85 20.3  /
  \multiput{\bulletsymb} at 20 19.7 /
  
  \multiput{$\smallsq3$} at 23.85 19.4  23.85 19.7  23.85 20  23.85 20.3
  23.85 20.6 / 
  
  \multiput{$\smallsq3$} at 27.85 19.4  27.85 19.7  27.85 20  27.85 20.3
  27.85 20.6  27.85 20.9 /
  \multiput{\bulletsymb} at 28 19.4 /
  
  \multiput{$\smallsq3$} at 31.85 19.7 /
  \multiput{\bulletsymb} at 32 19.7 /
  
  \multiput{$\smallsq3$} at 35.85 20 /
  
  \multiput{$\smallsq3$} at 39.85 19.4  39.85 19.7  39.85 20  39.85 20.3 /
  \multiput{\bulletsymb} at 40 20 /
  \multiput{$\smallsq3$} at 25.85 22.4  25.85 22.7  25.85 23  25.85 23.3
  25.85 23.6  25.85 23.9 /
  \arr{24.5 -1} {25.5 -3}
  \arr{28.5 -1} {29.5 -3}
  \arr{32.5 -1} {33.5 -3}
  \arr{26.5 -3} {27.5 -1}
  \arr{30.5 -3} {31.5 -1}
  \arr{34.5 -3} {35.5 -1}
  \arr{0.5 1} {1.25 2.5}
  \arr{4.5 1} {5.25 2.5}
  \arr{8.5 1} {9.25 2.5}
  \arr{12.5 1} {13.25 2.5}
  \arr{16.5 1} {17.25 2.5}
  \arr{20.5 1} {21.25 2.5}
  \arr{24.5 1} {25.25 2.5}
  \arr{28.5 1} {29.25 2.5}
  \arr{32.5 1} {33.25 2.5}
  \arr{36.5 1} {37.25 2.5}
  \arr{2.75 2.5} {3.5 1}
  \arr{6.75 2.5} {7.5 1}
  \arr{10.75 2.5} {11.5 1}
  \arr{14.75 2.5} {15.5 1}
  \arr{18.75 2.5} {19.5 1}
  \arr{22.75 2.5} {23.5 1}
  \arr{26.75 2.5} {27.5 1}
  \arr{30.75 2.5} {31.5 1}
  \arr{34.75 2.5} {35.5 1}
  \arr{38.75 2.5} {39.5 1}
  \arr{0.75 6.5} {1.5 5}
  \arr{4.75 6.5} {5.5 5}
  \arr{8.75 6.5} {9.5 5}
  \arr{12.75 6.5} {13.5 5}
  \arr{16.75 6.5} {17.5 5}
  \arr{20.75 6.5} {21.5 5}
  \arr{24.75 6.5} {25.5 5}
  \arr{28.75 6.5} {29.5 5}
  \arr{32.75 6.5} {33.5 5}
  \arr{36.75 6.5} {37.5 5}
  \arr{2.5 5} {3.25 6.5}
  \arr{6.5 5} {7.25 6.5}
  \arr{10.5 5} {11.25 6.5}
  \arr{14.5 5} {15.25 6.5}
  \arr{18.5 5} {19.25 6.5}
  \arr{22.5 5} {23.25 6.5}
  \arr{26.5 5} {27.25 6.5}
  \arr{30.5 5} {31.25 6.5}
  \arr{34.5 5} {35.25 6.5}
  \arr{38.6 5.2} {39.25 6.5}
  \arr{1 7.5} {1.5 7.25}
  \arr{5 7.5} {5.4 7.3}
  \arr{9 7.5} {9.5 7.25}
  \arr{13 7.5} {13.4 7.3}
  \arr{17 7.5} {17.5 7.25}
  \arr{21 7.5} {21.5 7.25}
  \arr{25 7.5} {25.5 7.25}
  \arr{29 7.5} {29.5 7.25}
  \arr{33 7.5} {33.5 7.25}
  \arr{37 7.5} {37.5 7.25}
  \arr{2.5 7.25} {3 7.5}
  \arr{6.5 7.25} {7 7.5}
  \arr{10.5 7.25} {11 7.5}
  \arr{14.5 7.25} {15 7.5}
  \arr{18.5 7.25} {19 7.5}
  \arr{22.5 7.25} {23 7.5}
  \arr{26.5 7.25} {27 7.5}
  \arr{30.5 7.25} {31 7.5}
  \arr{34.5 7.25} {35 7.5}
  \arr{38.5 7.25} {39 7.5}
  \arr{0.67 9} {1.33 10}
  \arr{4.67 9} {5.33 10}
  \arr{8.67 9} {9.33 10}
  \arr{12.67 9} {13.23 9.85}
  \arr{16.67 9} {17.33 10}
  \arr{20.67 9} {21.33 10}
  \arr{24.67 9} {25.33 10}
  \arr{28.67 9} {29.33 10}
  \arr{32.67 9} {33.33 10}
  \arr{36.67 9} {37.33 10}
  \arr{2.77 9.85} {3.33 9}
  \arr{6.67 10} {7.33 9}
  \arr{10.67 10} {11.33 9}
  \arr{14.77 9.85} {15.33 9}
  \arr{18.67 10} {19.33 9}
  \arr{22.67 10} {23.33 9}
  \arr{26.67 10} {27.33 9}
  \arr{30.67 10} {31.33 9}
  \arr{34.67 10} {35.23 9.15}
  \arr{38.67 10} {39.33 9}
  \arr{0.67 13} {1.33 12}
  \arr{4.67 13} {5.33 12}
  \arr{8.67 13} {9.33 12}
  \arr{12.67 13} {13.33 12}
  \arr{16.67 13} {17.33 12}
  \arr{20.67 13} {21.33 12}
  \arr{24.67 13} {25.33 12}
  \arr{28.67 13} {29.33 12}
  \arr{32.67 13} {33.33 12}
  \arr{36.67 13} {37.33 12}
  \arr{2.67 12} {3.33 13}
  \arr{6.67 12} {7.33 13}
  \arr{10.67 12} {11.33 13}
  \arr{14.67 12} {15.33 13}
  \arr{18.67 12} {19.33 13}
  \arr{22.67 12} {23.33 13}
  \arr{26.67 12} {27.33 13}
  \arr{30.67 12} {31.33 13}
  \arr{34.67 12} {35.33 13}
  \arr{38.67 12} {39.33 13}
  \arr{0.67 15} {1.33 16}
  \arr{4.67 15} {5.33 16}
  \arr{8.67 15} {9.33 16}
  \arr{12.67 15} {13.33 16}
  \arr{16.67 15} {17.33 16}
  \arr{20.67 15} {21.33 16}
  \arr{24.67 15} {25.33 16}
  \arr{28.67 15} {29.33 16}
  \arr{32.67 15} {33.33 16}
  \arr{36.67 15} {37.33 16}
  \arr{2.67 16} {3.33 15}
  \arr{6.67 16} {7.33 15}
  \arr{10.67 16} {11.33 15}
  \arr{14.67 16} {15.33 15}
  \arr{18.67 16} {19.33 15}
  \arr{22.67 16} {23.33 15}
  \arr{26.67 16} {27.33 15}
  \arr{30.67 16} {31.33 15}
  \arr{34.67 16} {35.33 15}
  \arr{38.67 16} {39.33 15}
  \arr{0.5 19.25} {1.33 18}
  \arr{4.5 19.25} {5.33 18}
  \arr{8.5 19.25} {9.33 18}
  \arr{12.5 19.25} {13.33 18}
  \arr{16.5 19.25} {17.33 18}
  \arr{20.5 19.25} {21.33 18}
  \arr{24.5 19.25} {25.33 18}
  \arr{28.5 19.25} {29.33 18}
  \arr{32.5 19.25} {33.33 18}
  \arr{36.5 19.25} {37.33 18}
  \arr{2.67 18} {3.5 19.25}
  \arr{6.67 18} {7.5 19.25}
  \arr{10.67 18} {11.5 19.25}
  \arr{14.67 18} {15.5 19.25}
  \arr{18.67 18} {19.5 19.25}
  \arr{22.67 18} {23.5 19.25}
  \arr{26.67 18} {27.5 19.25}
  \arr{30.67 18} {31.5 19.25}
  \arr{34.67 18} {35.5 19.25}
  \arr{38.67 18} {39.5 19.25}
  \arr{24.5 20.75} {25.5 22.25}
  \arr{26.5 22.25} {27.5 20.75}
  \setdots<2pt>
  \plot 0.5 -.15  3.5 -.15 /
  \plot 4.5 -.15  7.5 -.15 /
  \plot 8.5 -.15  11.5 -.15 /
  \plot 12.5 -.15  15.5 -.15 /
  \plot 16.5 -.15  19.5 -.15 /
  \plot 20.5 -.15  23.5 -.15 /
  \plot 26.5 -4.15  29.5 -4.15 /
  \plot 30.5 -4.15  33.5 -4.15 /
  \plot 36.5 -.15  39.5 -.15 /
  \plot 0.5 19.85  3.5 19.85 /
  \plot 4.5 19.85  7.5 19.85 /
  \plot 8.5 19.85  11.5 19.85 /
  \plot 12.5 19.85  15.5 19.85 /
  \plot 16.5 19.85  19.5 19.85 /
  \plot 20.5 19.85  23.5 19.85 /
  \plot 28.5 19.85  31.5 19.85 /
  \plot 32.5 19.85  35.5 19.85 /
  \plot 36.5 19.85  39.5 19.85 /
  
  \setsolid 
  \plot 0 1.2  0 6.6 /
  \plot 0 9  0 12.6 /
  \plot 0 15  0 18.6 /
  \plot 40 1.5  40 6.9 /
  \plot 40 9.3  40 12.9 /
  \plot 40 15.3  40 18.9 /
  \setshadegrid span <1.5mm>
  \vshade   0   -.15 19.85  <,z,,> 
  24  -.15 19.85  <z,z,,>
  26 -4    22.85  <z,z,,>
  28 -4    19.85  <z,z,,>
  34 -4    19.85  <z,z,,>
  36  -.15 19.85  <z,z,,>
  40  -.15 19.85 /
  \endpicture}
\begin{figure}[ht]  \label{arq36}
  \rotatebox{90}{\arqthreesix}
\end{figure}
\afterpage{\clearpage}

    Thus the matrix units or radical generators
  form the four generators of $\rad\End_{\mathcal S}X$ as a $\Lambda$-module.
  For each choice of $X$, each generator of $\rad\End_{\mathcal S}X$ factors through the
  middle term of the sequence, thus condition $2.$ or $2'.$ in Proposition~\ref{ar-test}
  is satisfied.  As a consequence, the sequences labelled $(2)$, $(3)$, $(4)$,
  and hence also the sequences ending at $\tau^2 C_2$, $\tau C_2$, $\tau^{-1}C_2$,
  $\tau C_3$ and $\tau^{-1}C_3$ 
  are Auslander-Reiten sequences for which the middle term decomposes as indicated. 

  \smallskip The module $A=(U\subset V)$  at the center of the slice 
  where the summands of $V=V_1\oplus V_2\oplus V_3$ 
  have length $6,4,2$,   is indecomposable:
  Consider the following submodules of $V$ which are invariant under automorphisms of $A$.
  By $\rho$ we denote the endomorphism of $V$ given by a radical generator of $\Lambda$.
  \begin{eqnarray*} L_1(U\subset V) & = & \rad^4\!V\cap \soc V \\
    L_2(U\subset V) & = & (\rad^3\!V\cap \soc^2\!V)\\
    L_i(U\subset V) & = & \rho^{-(i-2)}(L_2(U\subset V))\qquad \qquad{\rm for}\; i>2
  \end{eqnarray*}
  $$
  \beginpicture 
  \setcoordinatesystem units <.4cm,.4cm>
  \multiput{\sq} at 0 5  0 4  0 3  0 2  0 1  0 0  1 4  1 3  1 2  1 1  
  2 3  2 2  /
  \multiput{$\bullet$} at 0.3 1.8  1.3 1.8  1.7 2.2  2.7 2.2 /
  \plot 0.3 1.8  1.3 1.8 /
  \plot 1.7 2.2  2.7 2.2 /
  \setdashes <2pt>
  \plot -2 0.5  -.5 0.5 /
  \plot 2.5 0.5  6 0.5 /
  \plot -2 1.5  -.5 1.5 /
  \plot 3.5 1.5  6 1.5 /
  \plot -2 2.5  -.5 2.5 /
  \plot 3.5 2.5  6 2.5 /
  \plot -2 3.5  -.5 3.5 /
  \plot 3.5 3.5  6 3.5 /
  \plot -2 4.5  -.5 4.5 /
  \plot 2.5 4.5  6 4.5 /
  \put{$L_1$} at 6.5 0.5 
  \put{$L_2$} at 6.7 1.5 
  \put{$L_3$} at 6.5 2.5 
  \put{$L_4$} at 6.7 3.5 
  \put{$L_5$} at 6.5 4.5 
  \endpicture
  $$
  As an isomorphism invariant for $A$ we have the following 
  system   of vectorspaces and linear maps.
  If $k=\Lambda/\rad\Lambda$ 
  and $\frac{L_3}{L_2}=\frac{\rad^3\!V_1}{\rad^4\!V_1}\oplus\frac{\rad^2\!V_2}{\rad^3\!V_2}
  \oplus \frac{\rad V_3}{0}=k\oplus k\oplus k$ then the map 
  $u:\frac{(U\cap L_3)+L_2}{L_2}\to \frac{L_3}{L_2}$ given by the inclusion of the
  subspace $k(1,1,0)\oplus k(0,1,1)$ in $k\oplus k\oplus k$ makes the system,
  which has dimension type ${\phantom{12}2\phantom{21} \atop 12321}$, indecomposable.
  Thus, $A$ is indecomposable.
  $$\CD @. @. \frac{(U\cap L_3)+L_2}{L_2} \\
  @. @. @VV u V \\
  L_1 @<\rho<< \frac{L_2}{L_1} @<\rho<< \frac{L_3}{L_2} 
      @<\rho<< \frac{L_5}{L_4} @<\rho<< \frac{V}{L_5}
  \endCD$$

  \smallskip
  Finally consider the Auslander-Reiten sequence 
  $\mathcal E:0\to A\to B\to C\to 0$ starting at $A$. 
  Our study of the sequences ending and starting at $C_3$ shows that $C$ is as indicated
  and that $C_3$ occurs as a direct summand of $B$.  Since $A$ occurs as summand of the
  middle term of the sequence starting at $A_4$, it follows that $\tau^{-1} A_4$ occurs as
  another summand of $B$.  By Theorem~\ref{almost-split}, the sequence $\mathcal E$
  is split exact in each component, and hence $B$ must have a third direct summand, say $B'$,
  of type $(4;2)$.  It remains to note that the arrow $A\to B'$ has trivial valuation;
  for this observe that the types of the summands of the middle term 
  of the Auslander-Reiten sequence  ending at $A$ are as specified.
  
  \smallskip
  We have seen that the shaded region in the diagram forms a slice of the Auslander-Reiten
  quiver for $\mathcal S_3(\Lambda)$.  According to Theorem~\ref{introthm}, this 
  translation quiver is independent of the choice of the base ring $\Lambda$ which can be
  any commutative uniserial ring of length 6. 
\end{proof}

\bigskip\bigskip
Address of the author:\hfill
\begin{minipage}[t]{7cm}
  Mathematical Sciences\\
  Florida Atlantic University\\
  Boca Raton, Florida 33431-0991\\
  United States of America\\[1ex]
  e-mail: {\tt markus@math.fau.edu}
\end{minipage}


\smallskip
\begin{thebibliography}{99}

\bibitem{af} F.~W.~Anderson, K.~R.~Fuller:
  \textit{Rings and Categories of Modules,}
  Springer GTM {\bf 13}, Second edition 1992.

\bibitem{ars} M.~Auslander, I.~Reiten, S.O.~Smal\o:
  \textit{Representation Theory of Artin Algebras,}
  Cambridge studies in advanced mathematics {\bf 36}, 1995.

\bibitem{as} 
  M.~Auslander and S.~O.~Smal\o, 
  \textit{Almost split sequences in subcategories,}
  Journal of Algebra \textbf{69} (1981), 426-454.

\bibitem{avino-bautista}
  M.~Avi\~n\'o Dias and R.~Bautista Ramos,
  \textit{The additive structure of indecomposable
    $\mathbb Z_{p^n}C_p$-modules,}
  Communications in Algebra, \textbf{24} (1996), 2567-2595.

\bibitem{bhw}
  D.~Beers, R.~Hunter, and E.~Walker,
  \textit{Finite valuated p-groups,} 
  Abelian Group Theory, Springer LNM {\bf 1006} (1983), 471--507.

\bibitem{birkhoff}
  G.~Birkhoff, \textit{Subgroups of abelian groups},
  Proc.~Lond.~Math.~Soc., II. Ser.\ \textbf{38} (1934), 385--401.

\bibitem{bongartz-gabriel}
  K.~Bongartz and P.~Gabriel, 
  \textit{Covering spaces in representation theory,}
  Inv.~Math.\ {\bf 65} (1982), 331-378.

\bibitem{ps}
  C.~Petroro and M.~Schmidmeier, 
  \textit{Abelian groups with a $p^2$-bounded subgroup, revisited,}
  manuscript (2006), 15pp, arXiv: {\tt arxiv.org/abs/math/0605664}.

\bibitem{richman} 
  F.~Richman and E.~A.~Walker,  
  \textit{Subgroups of p${}^5$-bounded groups,}
  in: Abelian groups and modules, Trends Math., Birkh\"auser, Basel, 1999, 
  pp.~55--73.
  

\bibitem{rs-wild}
  C.~M.~Ringel and M.~Schmidmeier, 
  \textit{Submodule categories of wild  repre\-sentation type},
  J.~Pure Appl.~Algebra \textbf{205} (2006), 412-422.
  
\bibitem{rs-inv}
  C.~M.~Ringel and M.~Schmidmeier,
  \textit{Invariant subspaces of nilpotent linear operators. I,}
  manuscript (2006), 55pp, arXiv: {\tt arxiv.org/abs/math/0608666.}

\bibitem{rump} W.~Rump,
  \textit{Irreduzible und unzerlegbare Darstellungen 
    klassischer Ordnungen,} 
  Bayreuther Mathematische Schriften \textbf{32}, Bayreuth, 1990. 
  
  
\bibitem{ms-bounded}
  M.~Schmidmeier, \textit{Bounded submodules of modules,}
  J.~Pure Appl.~Algebra \textbf{203}, (2005), 45-82.
  
\bibitem{oberwolfach} M.~Schmidmeier, 
  \textit{A Remark by M.C.R. Butler on Subgroup Embeddings,}
  Oberwolfach Report \textbf6 (2005), 380-382.

\bibitem{simson} D.~Simson, \textit{Chain categories of modules and 
    subprojective representations 
    of posets over uniserial algebras,} 
  Rocky Mountain J.~Math.~\textbf{32} (2002), 1627--1650.  
  
\end{thebibliography}
\end{document}